\documentclass[12pt]{article}
\usepackage{amssymb}

\begin{document}
\newcommand{\bea}{\begin{eqnarray}}
\newcommand{\ena}{\end{eqnarray}}
\newcommand{\beas}{\begin{eqnarray*}}
\newcommand{\enas}{\end{eqnarray*}}
\newcommand{\beq}{\begin{equation}}
\newcommand{\enq}{\end{equation}}
\newcommand{\ww}{{\tilde a}}
\newcommand{\Opttwo}{\mathbf{2}}
\newcommand{\Optone}{\mathbf{1}}
\newcommand{\Vp}{V_n^p}
\newcommand{\ester}{\vartheta}
\newcommand{\B}{{\cal L}}
%\newcommand{\cstar}{d_\alpha}
%\newcommand{\astar}{b_\alpha}
%\newcommand{\delt}{\eta}
%\mathbf{1}\mathcal{1}\mathrm{1}
\newcommand{\X}{\chi}
\newcommand{\From}{From }
\newcommand{\U}{Z^+}
\newcommand{\s}{r}
\def\qed{\hfill \mbox{\rule{0.5em}{0.5em}}}
\newcommand{\ignore}[1]{}
\newtheorem{theorem}{Theorem}[section]
\newtheorem{corollary}{Corollary}[section]
\newtheorem{conjecture}{Conjecture}[section]
\newtheorem{proposition}{Proposition}[section]
\newtheorem{lemma}{Lemma}[section]
\newtheorem{definition}{Definition}[section]
\newtheorem{example}{Example}[section]
\newtheorem{remark}{Remark}[section]
\newtheorem{case}{Case}[section]

\title{{\bf\Large Two Choice Optimal Stopping}\,\,}
\author{David Assaf$\,^0$, Larry Goldstein and Ester
Samuel-Cahn$^0$\\Hebrew University, University of
Southern California and Hebrew University}

\maketitle

\date{}

\begin{abstract}

Let $X_n,\ldots,X_1$ be i.i.d. random variables with distribution
function $F$. A statistician, knowing $F$, observes the $X$ values
sequentially and is given two chances to choose $X$'s using
stopping rules. The statistician's goal is to stop at a value of
$X$ as small as possible. Let $V_n^2$ equal the expectation of the
smaller of the two values chosen by the statistician when
proceeding optimally. We obtain the asymptotic behavior of the
sequence $V_n^2$ for a large class of $F$'s belonging to the
domain of attraction (for the minimum) ${\cal D}(G^\alpha)$, where
$G^\alpha(x)=[1-\exp(-x^\alpha)]{\bf I}(x \ge 0)$. The results are
compared with those for the asymptotic behavior of the classical
one choice value sequence $V_n^1$, as well as with the ``prophet
value" sequence $\Vp=E(\min\{X_n,\ldots,X_1\})$.
\end{abstract}
\footnotetext{This research was supported by the ISRAEL SCIENCE
FOUNDATION (grant number 879/01)}
\section{Introduction}
\label{intro}
Kennedy and Kertz (1990, 1991) study the asymptotic
behavior of the value sequence, as $n \rightarrow \infty$, when
optimally stopping an $n$ long sequence of i.i.d. random variables
with common distribution function $F$, with the objective being to
stop on as large a value as possible. They show that the
asymptotic behavior of the value sequence depends upon the domain
of attraction, for the maximum, to which $F$ belongs.

Recently Assaf and Samuel-Cahn (2000) and Assaf, Goldstein, and Samuel Cahn (2002)
have studied optimal stopping problems where the statistician is given several choices,
and his return is the expected value of the maximal element chosen. The goals in these works
were the derivation of ``prophet inequalities."

In the present paper we study the limiting behavior of the value
sequence when the statistician, knowing $F$, is given two choices.
It turns out to be more convenient here to take as objective to
stop on as small a value as possible, and therefore to take as the
statistician's goal the minimization of the expected value upon
stopping. In particular, we consider a situation where the
statistician would like to choose the smallest possible value from
the $n$ i.i.d variables $X_n,\ldots,X_1$ presented sequentially,
and, with the luxury of two choices, can take a first choice as a
`fallback' value to use in case that none of the remaining
variables are small enough to take as a second choice.

The two choice problem we consider is more difficult by an order
of magnitude than the optimal one-choice problem. To be convinced
of this, let $V_n^1(x)$ (which we will also denote by $g_n(x)$)
and $V_n^2(x)$ be the value of the optimal one and two choice
policy respectively, when applied to the i.i.d. sequence $X_n,
\ldots, X_1$, when the statistician is already guaranteed the
value $x$. Note that for convenience we are indexing the
variables so that the first one observed is $X_n$ and the last is
$X_1$. Then by the dynamic programming principle, for one choice
$V_1^1(x)=E[X_1 \wedge x]$ and  we have \bea \label{1.7}
V^1_{n+1}(x)=E[X_{n+1} \wedge V^1_n(x)] \quad \mbox{for $n \ge
1$,} \ena whereas with two choices, $V_2^2(x) = E[X_2 \wedge X_1
\wedge x]$ and, for $n \ge 2$,
\bea
\label{1.9}
V_{n+1}^2(x) &=& E[V^1_n(X_{n+1}) \wedge V_n^2(x)].
\ena

The first term inside the square brackets in (\ref{1.9})
corresponds to choosing the current variable $X_{n+1}$ and being
left with only one additional choice among the remaining $n$
observations, while the second term corresponds to passing up the
current random variable $X_{n+1}$ and retaining two choices, with
the guaranteed bound $x$, among the remaining $n$ observations.

Comparing (\ref{1.7}) and (\ref{1.9}) we see that for one choice
the expectation computed in (\ref{1.7}) is with respect to the
random variables $X_{n+1}$ with identical distributions, whereas
the distribution of the random variable $V_n^1(X_{n+1})$ in
(\ref{1.9}) depends on the function $V_n^1$ which changes with $n$
even though the sequence $X_n,\ldots,X_1$ is identically
distributed.

Let \beq \label{defxF} x_F=\sup\{x:F(x)<1\}. \enq When nothing is
guaranteed, the value for the one and two stop problems will be
denoted $V_n^1$ and $V_n^2$ respectively, and satisfy
$V_n^1=V_n^1(x_F)$ and $V_n^2=V_n^2(x_F)$.

The optimal stopping rules can be specified in the one and two
stop cases by the values $V_n^1$, and the values $V_n^2$ and
functions $V_n^1(x)$, respectively, as follows. For the one stop
case, if $X_{n+1}$ is smaller than $V_n^1$ the variable $X_{n+1}$
should be taken. For the two stop case, if $V_n^1(X_{n+1})<V_n^2$
then the variable $X_{n+1}$ should be taken as the first choice,
and the optimal one stop strategy then used on the remaining $n$
variables when there is a guaranteed upper bound of $X_{n+1}$. In
other words, if $X_{m_1}$ has already been chosen as the first
choice, then take $X_m, m < m_1$ as the second choice when
$X_m<V_{m-1}^1(X_{m_1})$.

As in the one choice problem, the asymptotic behavior of the value sequence depends
on which of the three extreme value classes the distribution function $F$ belongs to.
In the present paper, we only consider $F$ which belongs to one of these domains of
attraction and take up the study of the remaining two classes in subsequent work.
Specifically in this paper, by a suitable shift of the origin, we assume that the
distribution function $F$ of the i.i.d. random variables belongs to the domain of
attraction (for the minimum)
${\cal D}(G^\alpha)$, where $\alpha>0$ and
\bea
\label{Galpha}
G^{\alpha}(x) = \left\{ \begin{array}{ccc}
0 & \mbox{$x<0$} \\
1-\exp(-x^\alpha) & \mbox{$x  \ge  0$,}
\end{array}\right.
\ena and satisfies $F(0)=0$ and $F(x)>0$ for all $x>0$. (This is
the Type III of Leadbetter, Lindgren and Rootz\'en, 1983, and Type
$\Psi_\alpha$ of Resnick, 1987.) A necessary and sufficient
condition for $F \in {\cal D}(G^\alpha)$ is \beas
%\label{FinDalpha}
F(x) = x^\alpha L(x), \quad \mbox{where $L(x)$ is slowly varying
at 0, i.e.} \quad \lim_{x \downarrow 0} \frac{L(tx)}{L(x)} = 1
\quad \mbox{for all $t>0$};
\enas
a sufficient (and close to necessary) condition is
\beas
%\label{1.6}
\lim_{x \downarrow 0} \frac{xF'(x)}{F(x)} = \alpha, \enas see
e.g. de Haan, 1976, Theorem 4.

Let $\Vp$ be the expected value of the minimum of $n$ i.i.d.
random variables. The results for the maximum (see e.g. Resnick
1987, Chapter 2.1) and the work of Kennedy and Kertz (1991)
translate for the minimum as follows: If $F \in {\cal
D}(G^\alpha)$, then
\bea \nonumber \lim_{n \rightarrow \infty} n F(\Vp) &=&
\Gamma(1+1/\alpha)^\alpha \quad \mbox{and}\\ \label{1.10} \lim_{n
\rightarrow \infty} n F(V_n^1) &=& (1+1/\alpha). \ena

Our main result for a statistician with two choices
is as follows.
\begin{theorem}
\label{nancy}
Let $X_n,\ldots,X_1$ be non-negative integrable i.i.d. random variables with distribution function
\bea
\label{wlog1}
F(x)=x^\alpha L(x) \quad \mbox{where $\lim_{x
\downarrow 0}L(x)$ exists and equals $\B \in (0,\infty)$.}
\ena
Then the optimal two choice value $V_n^2$ satisfies \beq
\label{1.12} \lim_{n \rightarrow \infty} n F(V_n^2) =
h^\alpha(b_\alpha) \enq where $b_\alpha>0$ is the unique solution
to \beq \label{1.13} \int_0^y h(u) du + (1/\alpha - y)h(y)=0, \enq
and $h(y)$ is the function \beq \label{1.11} h(y) = \left(
\frac{y}{1+\alpha y/(\alpha + 1)} \right)^{1/\alpha} \quad
\mbox{for $y \ge 0$}. \enq
\end{theorem}

The value $h(b_\alpha)$ depends only on $\alpha$ but
unfortunately, unlike the values (\ref{1.10}) cannot be given in
closed form in terms of $\alpha$.  A short table of the limiting
values (\ref{1.10}) and of $h^\alpha(b_\alpha)$ are given in Table
1. The performance improvement in having two choices over having
only one is substantial, in that the optimal stopping value
becomes much closer to that of the prophet. For example, for a
distribution with $\alpha=1$ such as the uniform, the limiting
values (for the minimum) for the statistician with one choice is
2, with two choices it is $1.165\ldots$, while the value for the
prophet is 1. More explicitly, with $n$ variables the optimal
value for a statistician with one choice is roughly $2/n$, for the
prophet it is roughly $1/n$, and for a statistician with two
choices it is $1.165 \ldots /n$.

The paper is organized as follows. In Section 2 we derive some
fundamental identities when $F$ belongs to the family \beq
\label{1.14} {\cal U}^\alpha(x) = \left\{
\begin{array}{cl}
0 & \mbox{for $x < 0$}\\
x^\alpha  & \mbox{for $0 \le x \le 1$}\\
1 & \mbox{for $x > 1$,}
\end{array} \right.
\enq for a fixed value of $\alpha > 0$; we also show heuristics
which explain the form of the function $h(y)$ of (\ref{1.11}). In
Section 3 we show that a particular sequence of functions $h_n$,
which determine $V_{n+1}^2$, converges to $h$. Section 4 contains
some general convergence results. In Section 5 we prove Theorem
\ref{nancy} for the special family (\ref{1.14}), and some results
concerning the finiteness of the limit of the moments of properly
scaled randomly selected values. In Section 6 Theorem \ref{nancy}
is generalized to a wide class of distributions in ${\cal
D}(G^\alpha)$. Section 7 contains numerical results presented in
Table 1, along with explanations and several additional remarks.
\section{The Fundamental Equations and Heuristics}
%\label{fundamental}
In general, for $X$ with distribution function
$F$, we let
$$
g(x) = E[X \wedge x].
$$
When $F(0)=0$, writing $g(x)=x-\int_0^x F(u)du$, we see easily
that $g(x)$ is positive and strictly increasing on the interval
$(0,x_F)$; hence the same is true for $g_{n+1}(x)=g(g_n(x))$.

In the remainder of this Section we consider $F={\cal U}^\alpha$
as in (\ref{1.14}), and in all the following we consider
$\alpha>0$ as fixed, to avoid the necessity of indexing quantities
by $\alpha$. For ${\cal U}^\alpha$ we have explicitly on the
interval $[0,1]$ \beq \label{2.1} g(x) = E[X \wedge x] = x -
\frac{x^{\alpha + 1}}{\alpha + 1}, \enq and with $g_1(x) = g(x)$,
\beq \label{2.2}
g_{n+1}(x)=g_n(x)-\frac{g_n(x)^{\alpha+1}}{\alpha+1}, \quad
\mbox{$n \ge 1$}. \enq Since a statistician with two choices does
at least as well as one with a single choice
$$
g_n(0) = 0 \le V_n^2 \le V_n^1 = g_n(1), \quad \mbox{$n \ge 2$}.
$$
As we are interested in the two choice case, we will henceforth
write $V_n$ to denote $V_n^2$ whenever convenient. Because the
function $g_n$ is strictly increasing on $[0,1]$, there exists a
unique number $b_n \in [0,1]$ satisfying \beq \label{2.3}
V_n=g_n(b_n). \enq We call $b_n$ the ``threshold value"  for the
following reason; by (\ref{1.9}) the statistician at stage $n+1$
will choose $X_{n+1}$ when $g_n(X_{n+1}) < V_n$, that is, when
$X_{n+1} < b_n$.

Since $b_n \in [0,1]$, $P(X>b_n)=1-b_n^\alpha$, and the basic
equation (\ref{1.9}) becomes \beq \label{2.4} V_{n+1} = \int_0
^{b_n} g_n(x)\alpha x^{\alpha - 1}dx +(1-b_n^\alpha)V_n, \quad
\mbox{$n \ge 2$.} \enq Letting $U_k$ be independent ${\cal
U}[0,1]$ variables, $U_k^{1/\alpha}$ has distribution ${\cal
U}^\alpha$, and hence we may begin recursion (\ref{2.4}) at
$V_2=E[U_2^{1/\alpha} \wedge U_1^{1/\alpha}]$. Scaling, \beq
\label{2.5} W_n=n^{1/\alpha}V_n, \quad B_n=n^{1/\alpha}b_n \enq
and \beq \label{2.6} g_n(x)=xf_n(n x^\alpha). \enq Since $g_n(x)$
is defined and positive for $0 < x \le 1$, the function $f_n(x)$
is defined and positive for $0 < x \le n$, and setting $f_n(0)=1$
makes $f_n(x)$ continuous as $x \downarrow 0$, since $g_n'(0)=1$.

Substituting (\ref{2.6}) into (\ref{2.4}) and making the change of
variable $y=n x^\alpha$ we obtain
$$
V_{n+1} = n^{-(1+1/\alpha)} \int_0^{B_n^\alpha}y^{1/\alpha}f_n(y)dy + (1-b_n^\alpha)V_n,
\quad \mbox{$n \ge 2$.}
$$
Multiplying by $n^{1/\alpha}$ and setting \beq \label{2.7} h_n(y)
= y^{1/\alpha}f_n(y) \enq we have \beq \label{2.8} \left(
\frac{n}{n+1} \right)^{1/\alpha} W_{n+1} = \frac{1}{n}\int_0
^{B_n^\alpha}h_n(y)dy + (1-b_n^\alpha)W_n, \quad \mbox{$n \ge 2$.}
\enq

By (\ref{2.3}),(\ref{2.5}),(\ref{2.6}) and (\ref{2.7}), \beq
\label{2.9} W_n = h_n(B_n^\alpha) \enq and we can now write
(\ref{2.8}) as our fundamental equation \beq \label{21star} W_m=c
\quad \mbox{and} \quad \left( \frac{n}{n+1} \right)^{1/\alpha}
W_{n+1} = \frac{1}{n}\int_0^n (h_n(y)\wedge W_n )\, dy \quad
\mbox{for $n \ge m$}, \enq with $m=2$ and $c=2^{1/\alpha}
E[U_2^{1/\alpha} \wedge U_1^{1/\alpha}]$. Later we allow for
arbitrary initial times $m \ge 1$ and any positive starting values
$c$.

The remainder of this Section is devoted to a heuristic argument
explaining (\ref{1.12}) and (\ref{1.13}), the appearance and form
of the function $h$ in (\ref{1.11}) and of Theorem \ref{nancy}.
Firstly, $((n+1)/n)^{1/\alpha} = 1 + 1/(\alpha n)+O(n^{-2})$, and
if $B_n^\alpha$ and the integral below remain bounded, we have
from (\ref{2.8})
$$
W_{n+1} - W_n = n^{-1} \int_0^{B_n^\alpha} h_n(y) dy + n^{-1}(1/\alpha - B_n^\alpha)W_n
+ O(n^{-2}).
$$
If $W_n \rightarrow d_\alpha$ such that
$W_n=d_\alpha+a/n+O(n^{-2})$, then $n (W_{n+1} - W_n) \to 0$, and
multiplying by $n$ we have \beq \label{2.10} 0 =
\int_0^{B_n^\alpha} h_n(y) dy + (1/\alpha - B_n^\alpha)W_n + o(1),
\enq and if $B_n^\alpha \rightarrow b_\alpha$ and $h_n \rightarrow
h$ as $n \rightarrow \infty$, then (\ref{2.10}) suggests
$$
%\label{2.11}
0 = \int_0^{b_\alpha}h(y) dy +(1/\alpha - b_\alpha)d_\alpha,
$$
where from (\ref{2.9}) also
$$
%\label{2.12}
d_\alpha=h(b_\alpha),
$$
which explains (\ref{1.12}) and (\ref{1.13}) of Theorem
\ref{nancy}.

By (\ref{2.7}), finding the limiting $h$ is equivalent to finding
the limiting $f$, since \beq \label{2.13} h(y) = y^{1/\alpha}f(y).
\enq Using (\ref{2.2}) and (\ref{2.6}) and the substitution $y=
nx^\alpha$, it follows that \beq \label{3.1}
f_{n+1}((1+\frac{1}{n})\,y) = f_n(y) - \frac{y}{(\alpha +
1)n}f_n(y)^{\alpha + 1}. \enq Subtracting $f_n(y)$ from both
sides, dividing by $y/n$ and taking limits as $n \rightarrow
\infty$ indicates that the limiting function $f$ should satisfy
the differential equation \beq \label{2.15} f'(y) = -f(y)^{\alpha
+ 1}/(\alpha + 1) \quad \mbox{with initial condition $f(0)=1$.}
\enq Equation (\ref{2.15}) has the unique solution \beq
\label{2.16} f(y) = (1 + \frac{\alpha y}{\alpha + 1})^{-1/\alpha},
\enq which together with (\ref{2.13}) yields the function $h$ of
(\ref{1.11}).

\section{Preliminary Lemmas}
%\label{lemmas}
In this Section we continue to consider $F={\cal U}^\alpha$ as in
(\ref{1.14}). With $f_n$ as in (\ref{2.6}) and $h_n$ as in
(\ref{2.7}), we have the following Lemma.
\begin{lemma}
%\label{lemma3.1}
The function $f_n(y)$ is strictly decreasing in $y$ for
$y \in [0,n]$ and $h_n(y)$ is strictly increasing in $y$ for $y \in [0,n]$.
\end{lemma}

\noindent {\bf Proof:}  We prove the lemma by induction on $n$.
For $n=1$ from (\ref{2.1}) and (\ref{2.6})
$$
f_1(y) = 1 - \frac{y}{\alpha + 1},
$$
so the result is immediate for $f_1$, and for $h_1$ by
(\ref{2.7}). Now assume the assertions are true for $n$. We shall show they are true
for $n+1$.  Note that for $0  \leq   y \leq  n$ we have
$0  \leq  y(n+1)/n  \leq  n+1$.  Differentiating (\ref{3.1}), for $0<y \le n$,
\beas
\left( \frac{n+1}{n} \right)f_{n+1}'(\frac{n+1}{n}y) &=&
f_n'(y) - \frac{1}{(\alpha + 1)n}f_n(y)^{\alpha + 1} -
\frac{y}{n} f_n^\alpha(y) f_n'(y)\\
&=& f_n'(y)[1-\frac{y}{n}f_n^\alpha (y)] - \frac{1}{(\alpha + 1)n}f_n(y) ^{\alpha + 1}\\
&<& f_n'(y)[1-\frac{y}{n}] - \frac{1}{(\alpha + 1)n}f_n(y) ^{\alpha + 1}\\
&<& 0,
\enas
where we have used
$f'_n(y)<0$ and $0 < f_n^\alpha(y) < 1$ for $0< y  \leq  n$.

\From (\ref{3.1}) and (\ref{2.7}) we have
$$
(\frac{n}{n+1})^{1/\alpha}h_{n+1}(\frac{n+1}{n}y)=
h_n(y)-\frac{1}{(\alpha+1)n}h_n(y)^{\alpha+1}.
$$
Thus for $0  \leq  y  <  n$ we have
$$
(\frac{n}{n+1})^{1/\alpha - 1}h_{n+1}'(\frac{n+1}{n}y)=
h_n'(y)[1-\frac{1}{n}h_n(y)^\alpha] > 0
$$
since by the induction hypothesis $h_n'(y) > 0$
and
\beas
h_n(y)^\alpha < h_n(n)^\alpha
 = [n^{1/\alpha}f_n(n)]^\alpha
 < n f_n^\alpha(0)
 = n. \hspace{1in} \qed
\enas
%$\qed$

Let $f(y)$ be given by (\ref{2.16}) and define
\beq
\label{3.3}
\epsilon_n(y) = f(y) - f_n(y).
\enq

\begin{lemma}
\label{lemma3.2}
With $\epsilon_n(y)$ as in (\ref{3.3}),
\beq
\label{3.4}
\epsilon_n(y) > 0 \quad \mbox{for $0 < y \le n.$}
\enq
\end{lemma}

\noindent {\bf Proof.}  We use the following two well known
inequalities. \beq \label{3.5} \mbox{For $0<\alpha \le 1$ and $x
\ge -1$,} \quad (1+x)^{\alpha} \le 1 + \alpha x, \enq and \beq
\label{3.6} \mbox{for $\alpha \ge 1$ and $x \ge -1$,} \quad 1 +
\alpha x \le (1+x)^{\alpha}. \enq We prove the lemma by
induction. For $n=1$ we must show that for $0 < y \le 1$
$$
1-\frac{y}{\alpha+1} < (1+\frac{\alpha y}{\alpha+1})^{-1/\alpha}
$$
which is equivalent to
$$
(1-\frac{y}{\alpha+1})^{\alpha} < (1+\frac{\alpha y}{\alpha+1})^{-1}
$$
or \beq \label{3.7} (1+\frac{\alpha
y}{\alpha+1})(1-\frac{y}{\alpha+1})^{\alpha} < 1. \enq Now for
$0 < \alpha  \leq  1$ we have by (\ref{3.5}) that the left hand
side of (\ref{3.7}) is less than or equal to
$$
(1+\frac{\alpha y}{\alpha+1})(1-\frac{\alpha y}{\alpha+1})=
1-(\frac{\alpha y}{\alpha+1})^2 < 1.
$$

For $\alpha>1$ the left hand side of (\ref{3.7}) is by (\ref{3.6}) less than
$$
(1+\frac{y}{\alpha+1})^{\alpha}(1-\frac{y}{\alpha+1})^{\alpha}=
[1-(\frac{y}{\alpha+1})^2]^{\alpha}<1.
$$
Thus $\epsilon_1(y) > 0$ for $0<y \le 1$.

Now suppose $\epsilon_n(y) > 0$ for $0<y \leq  n$. That
$\epsilon_{n+1}(y) > 0$ for $0<y \leq  n+1$, is equivalent to
$$
f_{n+1}(y) < (1+\frac{\alpha y}{\alpha+1})^{-1/\alpha}.
$$
By the induction hypothesis
$$
f_n(y) < (1+\frac{\alpha y}{\alpha+1})^{-1/\alpha} \quad
\mbox{for $0<y \leq  n$}
$$
and thus by (\ref{2.6})
$$
g_n(x) < x(1+\frac{\alpha n x^{\alpha}}{\alpha+1})^{-1/\alpha}
\quad \mbox{for $0 < x \le 1$},
$$
and since $g(\cdot)$ is an increasing function, using (\ref{2.2}), \beq \label{3.8}
g_{n+1}(x) < x(1+\frac{\alpha n x^\alpha}{\alpha +
1})^{-1/\alpha} [1-\frac{x^\alpha}{\alpha+1}(1+\frac{\alpha n
x^\alpha}{\alpha + 1})^{-1}]. \enq Thus, again by (\ref{2.6}), it
suffices to show that the right hand side of (\ref{3.8}) is less
than
$$
x(1+\frac{\alpha(n+1)x^{\alpha}}{\alpha+1})^{-1/\alpha}, \quad \mbox{for
$0<x \leq 1$.}
$$
Set $y=x^{\alpha}/(\alpha+1)$.  Then it suffices to show that
$$
(1+\alpha n y)^{-1/\alpha}[1-\frac{y}{1+\alpha n y}] <
(1+\alpha(n+1)y)^{-1/\alpha} \quad \mbox{for $0<y \le 1$,}
$$
i.e. that
$$
[1+\frac{\alpha y}{1+\alpha n
y}]^{1/\alpha}[1-\frac{y}{1+\alpha n y}] < 1,
$$
which is equivalent to \beq \label{3.9} [1+\frac{\alpha
y}{1+\alpha n y}][1-\frac{y}{1+\alpha n y}]^{\alpha} < 1. \enq

For $\alpha  \leq  1$ use (\ref{3.5}) to get that the left hand side of
(\ref{3.9}) is less than or equal to
$$
[1+\frac{\alpha y}{1+\alpha n y}] [1-\frac{\alpha y}{1+\alpha n
y}]=1-(\frac{\alpha y}{1+\alpha n y})^2<1.
$$
For $\alpha>1$ use (\ref{3.6}) to get that the left hand side of (\ref{3.9})
is less than
$$
[1+\frac{y}{1+\alpha n y}]^{\alpha} [1-\frac{y}{1+\alpha n
y}]^{\alpha}= [1-(\frac{y}{1+\alpha n y})^2]^\alpha <1. \qed
$$

\begin{lemma}
\label{lemma3.3}
With $\epsilon_n(y)$ as in (\ref{3.3}),
\beq
\label{3.10}
\epsilon_n(y) < \frac{y}{2n} \quad \mbox{for $0 < y \le n.$}
\enq
\end{lemma}

\noindent{\bf Proof:} We prove (\ref{3.10}) by induction. For
$n=1$ we must show that \beq \label{3.11} (1+\frac{\alpha
y}{\alpha+1})^{-1/\alpha} < 1-y(\frac{1}{\alpha+1}-\frac{1}{2})
\quad \mbox{for $0<y  \leq  1$.} \enq For $\alpha \ge 1$,
equation (\ref{3.11}) is obvious, since the left hand side is
less than 1 and the right hand side is greater than 1. For
$\alpha<1$ we have, by (\ref{3.6}) that
$$
(1+\frac{\alpha y}{\alpha+1})^{1/\alpha} \ge 1+\frac{y}{\alpha+1}.
$$
Thus to show (\ref{3.11}) it suffices to show
$$
\frac{1}{1+y/(\alpha+1)}<1-\frac{y(1-\alpha)}{2(\alpha+1)},
$$
i.e. that
$$
1<(1+\frac{y}{\alpha+1})
(1-\frac{y(1-\alpha)}{2(\alpha+1)})=1+\frac{y}{2}-\frac{y^2(1-\alpha)}{2(\alpha+1)^2}
$$
which clearly holds for $0<y \le 1$.

Now suppose (\ref{3.10}) holds for $n$. Let $0<y \le n+1$, and
$p_n=n/(n+1)$. By (\ref{3.1}) \beas
\epsilon_{n+1}(y) &=& f(y) - f_{n+1}(y) \\
&=& f(y) - f_n(p_n y) + \frac{p_n y}{(\alpha +1)n}f_n(p_n y)^{\alpha + 1}\\
&=& (f(y) - f(p_n y)) + (f(p_n y) - f_n(p_n y))  +
\frac{y}{(\alpha +1)(n+1)}f_n(p_n y)^{\alpha + 1}.
\enas
Thus
\beq
\label{3.12}
\epsilon_{n+1}(y) = f(y) - f(p_n y) + \epsilon_n(p_n y) + \frac{y}{(\alpha + 1)(n + 1)}
f_n(p_n y)^{\alpha + 1}.
\enq
Note that
\beq
\label{3.13}
f'(y)=-f(y)^{\alpha+1}/(\alpha+1)<0 \quad  \mbox{for $y>0$}
\enq
and
\beq
\label{3.14}
f''(y)=f(y)^{2\alpha+1}/(\alpha+1)>0 \quad \mbox{for $y>0$.}
\enq
Thus if we use the Taylor expansion
$$
f(x+\Delta) = f(x) +\Delta f'(x) + \frac{\Delta^2}{2}f''(x + \xi \Delta) \quad \mbox{for some $0 < \xi < 1$}
$$
with $x=p_ny$ and $\Delta=y/(n + 1)$ so that $x+\Delta=y$, we get,
by use of (\ref{3.13}) and (\ref{3.14}) \beq \label{3.15}
f(y)-f(p_ny) = -\frac{y}{(\alpha+1)(n+1)}f(p_ny)^{\alpha+1} +
\frac{y^2}{2(\alpha+1)(n+1)^2}f(\theta y)^{2\alpha+1} \enq where
$p_n<\theta<1$. Substituting (\ref{3.15}) into (\ref{3.12}) yields
\beq \label{3.16} \epsilon_{n+1}(y)=\epsilon_n(p_n
y)-\frac{y}{(\alpha+1)(n+1)}[f(p_ny)^{\alpha+1}-
f_n(p_ny)^{\alpha+1}]+\frac{y^2}{2(\alpha+1)(n+1)^2}f(\theta
y)^{2\alpha+1}. \enq Since by (\ref{3.4}) $f(p_ny)>f_n(p_ny)$ for
$0<y \leq  n+1$, we have \beq \label{3.17}
f(p_ny)^{\alpha+1}-f_n(p_n y)^{\alpha+1}>f(p_n
y)^{\alpha}[f(p_ny)-f_n(p_ny)]= f(p_n y)^{\alpha}\epsilon_n(p_ny).
\enq Substituting (\ref{3.17}) into (\ref{3.16}) yields \beq
\label{3.18} \epsilon_{n+1}(y) <  \epsilon_n (p_n
y)[1-\frac{y}{(\alpha+1)(n+1)} f(p_n
y)^{\alpha}]+\frac{y^2}{2(\alpha+1)(n+1)^2}f(\theta
y)^{2\alpha+1}. \enq It follows from the induction hypothesis that
for $0 <  y \leq  n+1$ (so that $0  <  p_n y  \leq  n$)
$$
\epsilon_n(p_n y)<\frac{p_n y}{2n}=\frac{y}{2(n+1)}.
$$
Thus (\ref{3.18}) yields
\beas
\epsilon_{n+1}(y)&<&\frac{y}{2(n+1)}
[1-\frac{y}{(\alpha+1)(n+1)}
f(p_n y)^{\alpha}]+\frac{y^2}{2(\alpha+1)(n+1)^2}f(\theta y)^{2\alpha+1}\\
&<&\frac{y}{2(n+1)}
[1-\frac{y f(p_ny)^\alpha\{1-f(\theta y)^{\alpha+1}\}}{(\alpha+1)(n+1)}
]<\frac{y}{2(n+1)},
\enas
where we have used the fact that $f$ is decreasing, $\theta>p_n$,
and $f<1$. $\qed$

\begin{corollary}
\label{corollary3.1}
\beas
f_n(y) \rightarrow f(y)
&=& (1+\frac{\alpha y}{\alpha + 1})^{-1/\alpha}
\quad \mbox{for all $y>0$, as $n \rightarrow \infty$}\\
h_n(y) \rightarrow h(y)
&=& \left(\frac{y}{1+\alpha y/(\alpha + 1)}\right)^{1/\alpha}
\quad \mbox{for all $y>0$, as $n \rightarrow \infty$.}
\enas
\end{corollary}

\begin{remark}
%\label{remark3.1}
Note that by (\ref{2.6}),(\ref{2.7}) and (\ref{1.7})
$$
h_n(n)=n^{1/\alpha}g_n(1)=n^{1/\alpha}V_n^1
$$
and thus, by (\ref{1.10})
$$
\lim_{n \rightarrow \infty}h_n(n)=[1+ 1/\alpha]^{1/\alpha}.
$$
On the other hand, we also have
$$
\lim_{y \rightarrow \infty}h(y) = [1+ 1/\alpha]^{1/\alpha}
$$
Thus, the convergence to $h$ in Corollary \ref{corollary3.1} satisfies
$$
\lim_{n \rightarrow \infty}h_n(n) =
\lim_{y \rightarrow \infty} \lim_{n \rightarrow \infty} h_n(y).
$$
\end{remark}

\section{Convergence of Recursions}
%\label{general}
To prove convergence of the sequence $W_n$ determined by the
recursion (\ref{21star}), we first study the behavior of a
sequence $Z_n$, whose values are given by the simpler recursion
\beq \label{n4.5} Z_m=c \quad \mbox{and} \quad \left(\frac{n}{n+1}
\right)^{1/\alpha} Z_{n+1} = \frac{1}{n} \int_0^n (q(y) \wedge
Z_n)dy \quad \mbox{for $n\ge m$}, \enq where the function in the
integral does not depend on $n$.

For $\alpha>0$ a fixed value and $q(\cdot)$ a given function,
define \beq \label{n4.1} Q(y) = \int_0^y q(u) du + (1/\alpha - y)
q(y). \enq
We prove the convergence of $Z_n$ under the following conditions:\\
(i) $q(0)=0$\\
(ii) $q(u)$ for $0 <u < \infty$  is non-decreasing everywhere and
strictly increasing and differentiable for $0<u<A$ where
$1/\alpha < A \le \infty$.\\
 (iii) There exists a
unique positive root $b \in (1/\alpha,A)$ to the equation
$Q(y)=0$.

\begin{lemma}
\label{nlemma4.2} Under conditions (i) and (ii), the function
$Q(\cdot)$ is strictly increasing for $0 < y < 1/\alpha$, strictly
decreasing for $1/\alpha < y < A$, and non-increasing for $A<y$.
Hence $Q(A)=\lim_{y \uparrow A}Q(y)$ exists and is in
$[-\infty,\infty)$, even when $A=\infty$, and (iii) holds if $Q(A)
< 0$.
\end{lemma}
\noindent {\bf Proof:} For $0 \le y_1<y_2<1/\alpha$
straightforward calculations yield
$$
Q(y_2)-Q(y_1) \ge (q(y_2)-q(y_1))(1/\alpha - y_2),
$$
and for $1/\alpha < y_1<y_2$,
$$
Q(y_2)-Q(y_1) \le (q(y_2)-q(y_1))(1/\alpha-y_1).
$$
The claims now follow directly. $\qed$ \\[1ex]

The main result of this Section is
\begin{theorem}
\label{theorem4.1} Let (i), (ii) and (iii) hold, and let $Z_n$ be
given by (\ref{n4.5}) with $m \ge 1$ any integer and $c \in
(0,\infty)$ any constant. Then the limit of $Z_n$ exists and \beas
\lim_{n \rightarrow \infty}Z_n = d, \enas where $d=q(b)$ with $b$
the unique root of $Q(y)=0$.
\end{theorem}

Lemma \ref{for-inf} is the crux of of the proof of Theorem
\ref{theorem4.1}.
\begin{lemma}
\label{for-inf} Assume that (i), (ii) and (iii) hold. Let $m \ge
1$ be any integer and $c \in (0,\infty)$ any constant, and
suppose that $Z_n$ for $n \ge m$ is defined by (\ref{n4.5}). Then
for every $\delta \in (0,\min\{q(A)-d,d-q(1/\alpha)\})$ there
there exists $\Delta>0$ and $n_0$ such that for all $n \ge n_0$,
\bea \label{smallz} \mbox{if $Z_n < d-\delta$ then $Z_{n+1} \ge
(1+\Delta/n)Z_n$,} \ena \bea \label{bigz} \mbox{if $Z_n>d+\delta$
then $Z_{n+1} \le (1-\Delta/n) Z_n$,} \ena \bea \label{stays<d}
\mbox{if $Z_n < d$ then $Z_{n+1} < d$, and} \ena \bea
\label{withindelta} \mbox{if $|Z_n-d| \le \delta$ then
$|Z_{n+1}-d| \le \delta$.} \quad \ena
\end{lemma}

\noindent {\bf Proof:} We have
$$
(1 + \frac{1}{n})^{1/\alpha} = 1 + \frac{1}{\alpha n} +
\frac{1}{\alpha}(\frac{1}{\alpha}-1)\frac{1}{2n^2} +
O_\alpha(n^{-3}),
$$
and hence for $\gamma>0$
\beq \label{withgamma}
(\frac{n+1}{n})^{1/\alpha} (1-\frac{1}{n \gamma}) =  1 -
\frac{1}{n}(\frac{1}{\gamma}-\frac{1}{\alpha})+\frac{1}{n^2}\left(
\frac{1}{2\alpha}(\frac{1}{\alpha}-1)-\frac{1}{\alpha \gamma
}\right)+O_{\alpha,\gamma}(n^{-3}),\enq where we write
$O_\lambda(f_n)$ to indicate a sequence bounded in absolute value
by $f_n$ times a constant depending only on $\lambda$, a
collection of parameters.

Define
$$
%\label{1star}
M(t) = \int_0 ^{q^{-1}(t)} \left(1-\frac{q(y)}{t} \right) dy \quad
\mbox{for $0 \le t<q(A)$.}
$$
From (\ref{n4.1}), $Q(b)=0$ and $d=q(b)$, we have
$$
%\label{2star}
M(d) = 1/\alpha.
$$
It is not hard to see that $M(t)$ is strictly increasing over its
range. Hence, setting $\Delta_1= (1/\alpha-M(d-\delta))/2$ and
$\Delta_2= (M(d+\delta)-1/\alpha)/2$  we have
$\Delta=\min\{\Delta_1,\Delta_2\}>0$. Now consider the function
$$
%\label{3star}
r_n(t) = \frac{1}{n}\int_0^n \left(\frac{q(y)}{t} \wedge 1
\right) dy = 1 - \frac{1}{n}\int_0^{q^{-1}(t) \wedge n} \left(1 -
\frac{q(y)}{t} \right) dy.
$$
Since $Z_m>0$ we have $Z_n>0$ for all $n \ge m$, and now by
(\ref{n4.5}) we have \bea \label{4star} Z_{n+1}/Z_n = \left(
\frac{n+1}{n} \right)^{1/\alpha} r_n(Z_n). \ena

By definition\beas r_n(t) = 1 - \frac{1}{n}M(t) \quad \mbox{for $0
\le t < q(n)$.} \enas

To prove (\ref{smallz}), assume $Z_n<d-\delta$. Since $r_n$ is
decreasing, using (\ref{4star}) and (\ref{withgamma}), we have
for all $n > q^{-1}(d-\delta)$, \beas Z_{n+1} &\ge&
Z_n (\frac{n+1}{n})^{1/\alpha}r_n(d - \delta) \\
&=& Z_n (\frac{n+1}{n})^{1/\alpha} (1 -
\frac{1}{n}M(d-\delta))\\
&=&(1+\frac{1}{n}(\frac{1}{\alpha}-M(d-\delta))+O_{\alpha,d-\delta}(n^{-2}))Z_n\\
&\ge& (1+\frac{\Delta_1}{n}) Z_n \ge (1+\frac{\Delta}{n}) Z_n
\enas for all $n$ sufficiently large, showing (\ref{smallz}).

Next we prove (\ref{bigz}). When $Z_n \ge d+\delta$, we have
similarly that for $n > q^{-1}(d + \delta)$, \beas Z_{n+1} &\le&
Z_n (\frac{n+1}{n})^{1/\alpha}r_n(d + \delta) \\
&=& Z_n (\frac{n+1}{n})^{1/\alpha} (1 -
\frac{1}{n}M(d+\delta))\\
&=&(1-\frac{1}{n}(M(d+\delta)-\frac{1}{\alpha})+O_{\alpha,d+\delta}(n^{-2}))Z_n\\
&\le& (1-\frac{\Delta_2}{n})Z_n \le (1-\frac{\Delta}{n})Z_n\enas
for all $n$ sufficiently large.

Turning now to (\ref{stays<d}) and (\ref{withindelta}), for $Z_n
\le d+\delta$, since $d+\delta <q(A)$, $\beta_n$ is well defined
by
$$
q(\beta_n)=Z_n.
$$
Now by (\ref{n4.5}) and (\ref{n4.1})
$$
\left( \frac{n}{n+1} \right)^{1/\alpha}Z_{n+1} =
\frac{1}{n}\left( \int_0^{\beta_n} q(y) dy +
(n-\beta_n)q(\beta_n)\right) = \frac{1}{n}Q(\beta_n) + (1 -
\frac{1}{\alpha n})Z_n;
$$
thus \beq \label{n4.9} Z_{n+1} = \left( 1 + \frac{1}{n}
\right)^{1/\alpha} \frac{1}{n} Q(\beta_n) + R_n Z_n \enq where
\beq \label{ester-again-insists} R_n=(1+\frac{1}{n})^{1/\alpha}(1
- \frac{1}{\alpha n}). \enq

Consider
$$ Q(q^{-1}(u)) = \int_0^{q^{-1}(u)}q(y)dy + (1/\alpha -
q^{-1}(u))u.
$$
Since $q^{-1}(u)$ is differentiable for $0 < u < q(A)$,
$$
\frac{d}{du}Q(q^{-1}(u)) = 1/\alpha - q^{-1}(u).
$$
Hence, evaluating $Q(q^{-1}(u))$ by a Taylor expansion around
$d$, and using $Q(b)=Q(q^{-1}(d))=0$, we obtain that there exists
some $\xi_{Z_n}$ between $d$ and $Z_n$ such that \beq
\label{n4.10} Q(\beta_n) = Q(q^{-1}(Z_n))=(Z_n-d)(1/\alpha -
q^{-1}(\xi_{Z_n})). \enq Subtracting $d$ from both sides of
(\ref{n4.9}) and using (\ref{n4.10}) we obtain \beq \label{n4.11}
Z_{n+1} - d = \left\{
1-(1+\frac{1}{n})^{1/\alpha}\frac{1}{n}(q^{-1}(\xi_{Z_n}) -
\frac{1}{\alpha}) \right\}(Z_n-d) + [R_n-1]Z_n. \enq

Take $n_1$ such that for all $n \ge n_1$
$$
(1+\frac{1}{n})^{1/\alpha} \frac{1}{n} (q^{-1}(d)-1/\alpha)) < 1.
$$
Then for $Z_n<d$ we have $\xi_{Z_n}<d$ and hence
$q^{-1}(\xi_{Z_n})<q^{-1}(d)$, and so
$$
0 < \left\{
1-(1+\frac{1}{n})^{1/\alpha}\frac{1}{n}(q^{-1}(\xi_{Z_n}) -
\frac{1}{\alpha}) \right\}.
$$
Hence the first term on the right hand side of (\ref{n4.11}) is
strictly negative. Next, there exists $n_2 \ge n_1$ so that for $n
\ge n_2$ we have $0<R_n<1$, by (\ref{ester-again-insists}) and
(\ref{withgamma}) with $\gamma=\alpha$. For such $n$ the second
term on the right hand side is also negative, and the sum of
these two terms is therefore negative. This proves
(\ref{stays<d}).

To consider (\ref{withindelta}) suppose that $|Z_n-d| \le \delta$.
Then $ |\xi_{Z_n} -d| \le \delta$, and therefore
$$
q^{-1}(d-\delta) \le q^{-1}(\xi_{Z_n}) \le q^{-1}(d+\delta).
$$
Hence, for all $n$ sufficiently large so that
$$
(1+\frac{1}{n})^{1/\alpha}
\frac{1}{n}\left(q^{-1}(d+\delta)-1/\alpha \right) \le 1,
$$
letting $\Delta_3=q^{-1}(d-\delta)-1/\alpha$, which is strictly
positive by choice of $\delta<d-q(1/\alpha)$, we have
$q^{-1}(\xi_{Z_n})-1/\alpha
\ge
\Delta_3$ and therefore
\bea
\label{again-insists} 0 \le \left\{
1-(1+\frac{1}{n})^{1/\alpha}\frac{1}{n}(q^{-1}(\xi_{Z_n}) -
\frac{1}{\alpha}) \right\} &\le& 1-\frac{\Delta_3}{n.} \ena
Further, from (\ref{ester-again-insists}), again using
(\ref{withgamma}) with $\gamma=\alpha$, there exists $K_\alpha$
such that
$$
|R_n-1| \le \frac{K_\alpha}{n^2}.
$$
Then for all $n$ so large that
$$
\frac{K_\alpha}{n}(d+\delta) \le \Delta_3 \delta
$$
we have, using (\ref{n4.11}) and (\ref{again-insists}),

\beas |Z_{n+1}-d| & \le & (1-\frac{\Delta_3}{n})|Z_n-d| + |R_n-1|Z_n\\
&\le& (1-\frac{\Delta_3}{n})\delta + \frac{K_\alpha}{n^2}(d+\delta)\\
&\le& \delta. \enas This proves (\ref{withindelta}). $\qed$
 \\[1ex]

\noindent {\bf Proof of Theorem \ref{theorem4.1}:} For $\delta
\in (0,\min\{q(A)-d,d-q(1/\alpha)\})$, let $\Delta$ and $n_0$ be
as in Lemma \ref{for-inf}. \\
Case I: $Z_{n_0} > d+\delta$. If $Z_n>d+\delta$ for all $n \ge
n_0$ then by (\ref{bigz}) we would have
$$
Z_{n+1} \le \prod_{j=n_0}^n(1-\frac{\Delta}{j})Z_{n_0} \rightarrow
0,
$$
a contradiction. Hence for some $n_1 \ge n_0$ we have $Z_{n_1} \le
d+\delta$, and we would therefore be in Case II or Case III.\\

Case II: $Z_{n_1} < d - \delta$ for some $n_1 \ge n_0$. If $Z_n <
d -\delta$ for all $n \ge n_1$ we would have by (\ref{smallz})
that
$$
Z_{n+1} \ge \prod_{j=n_1}^n(1+\frac{\Delta}{j})Z_{n_1} \rightarrow
\infty,
$$ a contradiction. Hence there exists $n_2 \ge n_1$ such that $Z_{n_2}\ge d - \delta$. By
(\ref{stays<d}), $Z_{n_2}<d$, reducing to Case III.\\

Case III: $|Z_{n_1}-d| \le \delta$ for some $n_1 \ge n_0$. In this
case $|Z_n-d| \le
\delta$ for all $n \ge n_1$, by (\ref{withindelta}).\\
Since $\delta$ can be taken arbitrarily small, the Theorem is
complete. $\qed$

The following Lemma may be of general interest, and presumably
has been noticed independently by others. We will apply it to
obtain asymptotic properties of moments in Section \ref{five}.
\begin{lemma}
\label{lemma4.1} \noindent A. Let $D_n, n \ge n_0$ be a
non-negative sequence satisfying \beq \label{ester-insists}
D_{n+1} \le \ester_n D_n + \gamma_n, \quad n \ge n_0, \enq where
$$
0 \le \ester_n \le (1-\ester/n) \quad \mbox{and} \quad 0 \le
\gamma_n \le \frac{C}{n}
$$
for some $\ester>0$ and $C \ge 0$. Then
$$
\limsup_{n \rightarrow \infty}D_n < \infty.
$$

\noindent B. Let $D_{n_0}>0$ and let  \beq
\label{didn't-insist} D_{n+1} \ge \ester_n D_n + \gamma_n, \quad
n \ge n_0, \enq where
$$
\ester_n \ge (1+\ester/n), \quad \mbox{and} \quad \gamma_n \ge 0.
$$
for some $\ester>0$.Then
$$
\lim_{n \rightarrow \infty} D_n = \infty.
$$

\end{lemma}
\noindent {\bf Proof:} Consider A. If (\ref{ester-insists})
holds, then by induction, for all $n \ge n_0$ and $k \ge 0$, \beq
\label{Dninduct} D_{n+k+1} \le \left( \prod_{j=n}^{n+k} \ester_j
\right) D_n + \sum_{j=n}^{n+k}\left( \prod_{l=j+1}^{n+k} \ester_l
\right) \gamma_j. \enq Using $\ester_n \le (1-\ester/n)$ and $1-x
\le e^{-x}$ we have \beas
\prod_{l=j+1}^{n+k} \vartheta_l &\le& \prod_{l=j+1}^{n+k} e^{-\ester/l} \\
                                &=& \exp(-\ester \sum_{l=j+1}^{n+k} 1/l )\\
                                &\le& \exp( -\ester (\log(n+k) - \log(j+1)))\\
                                &=& \left( \frac{j+1}{n+k} \right)^\ester.
\enas Hence, from (\ref{Dninduct}), for all $k \ge 0$, \beas
D_{n+k+1} &\le& \left(\prod_{j=n}^{n+k}\vartheta_j \right)D_n +
 \sum_{j=n}^{n+k} \left( \prod_{l=j+1}^{n+k}\vartheta_l \right) \gamma_j \\
 &\le& D_n + \sum_{j=n}^{n+k} \left( \frac{j+1}{n+k} \right)^\ester \frac{C}{j}\\
 &\le& D_n + \frac{2^\ester C}{(n+k)^\ester} \sum_{j=n}^{n+k} j^{\ester-1} \\
 &\le& D_n + \frac{2 ^\ester C}{\ester}\left(\frac{n+k+1}{n+k} \right)^\ester.
\enas Letting $k \rightarrow \infty$ we see that the $D_n$
sequence is bounded.

To prove B, note that $D_n>0$ for all $n \ge n_0$ and that for all
$j$ sufficiently large
$$
\ester_j \ge (1+\ester/j) \ge \exp( \ester/(2j)),
$$
which gives, by (\ref{didn't-insist}), \beas D_{n+k+1} \ge \left(
\prod_{j=n}^{n+k} \ester_j \right) D_n
          \ge \exp(\frac{\ester}{2} \sum_{j=n}^{n+k}\frac{1}{j})D_n \rightarrow \infty
          \quad \mbox{as $k \rightarrow \infty$}. \,\,\,\qed
\enas

\section{The Family ${\cal U}^\alpha$}
\label{five} As in (\ref{n4.1}), with $h(\cdot)$ defined in
(\ref{1.11}), let \beas H(y)=\int_0^y h(u) du + (1/\alpha - y)
h(y); \enas note that $h(\cdot)$ is strictly increasing for $0
\le y < \infty$.
\begin{lemma}
\label{nlemma5.1} There exists a unique value $b_\alpha >
1/\alpha$ such that $H(b_\alpha)=0$, and \beq \label{2<1}
h^\alpha(b_\alpha)<1+\frac{1}{\alpha}. \enq
\end{lemma}
\noindent {\bf Proof:} By Lemma \ref{nlemma4.2}, $H(y)$ is strictly increasing
for $0 < y < 1/\alpha$ and strictly decreasing for $1/\alpha < y < \infty$.
Hence a root exists in $(1/\alpha,\infty)$ and is unique if $H$ is ever
negative. Since
$$
H'(y) = (1/\alpha - y)h'(y),
$$
for some constant $a$
\beq
\label{Hasint}
H(y) = a + \int_{1/\alpha}^y (1/\alpha - u)h'(u)du.
\enq
Now, since $h(y)$ converges to a finite positive limit at infinity, and
$$
h'(y) = \frac{1}{\alpha}h(y)^{1-\alpha}\frac{1}{(1+\alpha y/(\alpha +1))^2},
$$
we have that $y^2 h'(y)$ is bounded away from zero and infinity as $y \rightarrow \infty$,
and therefore
$$
\int_{1/\alpha}^\infty h'(u)du < \infty \quad \mbox{and} \quad
\int_{1/\alpha}^y uh'(u) du \rightarrow \infty \quad \mbox{as $y \rightarrow \infty$,}
$$
yielding from (\ref{Hasint}) that
$$
\lim_{y \rightarrow \infty}H(y) = -\infty.
$$
Inequality (\ref{2<1}) follows from $\lim_{y \rightarrow
\infty}h^\alpha(y)= 1+1/\alpha$ .$\qed$ \\[1ex]

For $f(y)$ as given in (\ref{2.16}), setting \beq \label{n5.1}
f_j^*(y)=f(y) - y/2j \enq we have \beq \label{n5.2}
\frac{d}{dy}\left(y^{1/\alpha}f_j^*(y) \right) = y^{1/\alpha -
1}\left( \frac{f(y)}{\alpha}- \frac{yf(y)^{\alpha+1}}{\alpha +1}-
\frac{y}{2j}(1/\alpha + 1)\right). \enq Since $yf(y)^\alpha$ is
strictly increasing with limit $(\alpha+1)/\alpha$ at infinity, $
f(y)/\alpha > yf(y)^{\alpha+1}/(\alpha+1)$ for all $y \ge 0$.
Hence, for any fixed $A > b_\alpha$ we have
$$
\inf_{0 \le y \le A} \left( \frac{f(y)}{\alpha} -
\frac{yf(y)^{\alpha+1}}{\alpha+1} \right)>0.
$$ It follows that there exists $j_0=j_0(A)$ such that the
derivative in (\ref{n5.2}) is positive for all $0<y\le A$ and all
$j > j_0$. For these $j$, set \beq \label{n5.3} k_j(y) = \left\{
\begin{array}{ccc}
y^{1/\alpha} f_j^*(y) & \mbox{for $0 \le y < A$}\\
A^{1/\alpha} f_j^*(A) & \mbox{for $A \le y < \infty$}
\end{array} \right.
\enq
and
\beas
%\label{n5.5}
K_j(y) = \int_0^y k_j(u) du + (1/\alpha - y) k_j(y).
\enas

\begin{lemma}
\label{nlemma5.1*} There exists $j_1$ such that for all $j > j_1$
there are unique roots $b_{j,\alpha}$ to $K_j(y)=0$ and
$b_{j,\alpha}> 1/\alpha$. Setting $d_{j,\alpha}=k_j(b_{j,\alpha})$
we have \beq \label{n5.6} b_{j,\alpha} \rightarrow b_\alpha \quad
\mbox{and} \quad d_{j,\alpha} \rightarrow d_\alpha \quad \mbox{as
$j \rightarrow \infty$, where $d_\alpha = h(b_\alpha)$ .} \enq
\end{lemma}
\noindent {\bf Proof:} We apply Lemma \ref{nlemma4.2}. The
functions $k_j(\cdot)$ satisfy $k_j(0)=0$, are non-decreasing
everywhere and are strictly increasing and differentiable for $0
< y < A$. Further, $k_j(y)$ converges uniformly to $h(y)$ in
$[0,A]$, yielding the uniform convergence of $K_j(y)$ to $H(y)$
in $[0,A]$. Since $H$ is strictly decreasing in $(1/\alpha,
\infty)$, it follows that $H(A) < H(b_\alpha)=0$. Hence, since
$K_j(A) \rightarrow H(A)$ as $j \rightarrow \infty$, for all $j$
sufficiently large $K_j(A) < 0$. For such $j$ Lemma
\ref{nlemma4.2} now yields the existence of a unique root
$b_{j,\alpha} > 1/\alpha$ satisfying $K_j(b_{j,\alpha}) =0$.

The uniform convergence of $K_j$ to $H$ implies $H(b_{j,\alpha})
\rightarrow 0$ as $j \rightarrow \infty$, from which the
convergence of $b_{j,\alpha}$ to $b_\alpha$ follows. That
$d_{j,\alpha}$ converges to $d_\alpha$ follows from the uniform
convergence of $k_j$ to $h$ in $[0,A]$. $\qed$ \\[1ex]

\begin{lemma}
\label{nlemma5.2} Let $m \ge 1$ be any integer and $c \in
(0,\infty)$ be any constant. For $n > m$ let $W_n$ be determined
by the recursion (\ref{21star}) with starting value $W_m=c$, and
let \beq \label{11star} \U_m=c \quad \mbox{and} \quad
\left(\frac{n}{n+1} \right)^{1/\alpha} \U_{n+1}= \frac{1}{n}
\int_0^n \left(h(y) \wedge \U_n \right) dy \quad \mbox{for $n \ge
m$.} \enq With $j_1$ as in Lemma \ref{nlemma5.1*}, for all $j >
j_1$ let $m_j^*=\max\{m,j\}$. Now define sequences $Z_{j,n}^-$ for
$n \ge m_j^*$, by \beq \label{14star} Z_{j,m_j^*}^-=W_{m_j^*}
\quad \mbox{and} \quad \left( \frac{n}{n+1}
\right)^{1/\alpha}Z_{j,n+1}^- = \frac{1}{n}\int_0^n \left(k_j(y)
\wedge Z_{j,n}^- \right) dy \quad \mbox{for $n \ge m_j^*.$} \enq
Then for all $n \ge m_j^*$, \beq \label{12star} Z_{j,n}^- \le W_n
\le \U_n \enq and \beq \label{15star} \lim_{n \rightarrow \infty}
Z_{j,n}^- = d_{j,\alpha} \quad \mbox{and} \quad \lim_{n
\rightarrow \infty} Z_n^+ = d_\alpha. \enq
\end{lemma}

\noindent {\bf Proof:} With $j>j_1$ and $f_j^*$ defined in
(\ref{n5.1}), Lemmas \ref{lemma3.3}, \ref{lemma3.2} and
monotonicity of $f_n$ give
$$
f_j^*(y) < f_n(y) < f(y) \quad \mbox{for all $n \ge j$ and $0 < y
\le n$.}
$$
Therefore, by (\ref{n5.3}), (\ref{2.7}) and (\ref{2.13}),
$$
k_j(y) < h_n(y) < h(y) \quad \mbox{for all $n \ge j$ and $0 < y
\le n$.}
$$

Equation (\ref{12star}) now follows by a comparison of
(\ref{14star}), (\ref{21star}) and (\ref{11star}), and
(\ref{15star}) follows directly from Theorem \ref{theorem4.1}.
$\qed$

It is convenient to consider the value and scaled value arising
from stopping a sequence
$U^{1/\alpha}_n,\ldots,U^{1/\alpha}_{m+1},X_m,X_{m-1},\ldots,X_1$
of independent variables with a finite initial subsequence from a
distribution other than that of $U^{1/\alpha}$. The scaled value
sequence for this problem satisfies (\ref{21star}) with
$c=m^{1/\alpha}V_m(X_m,\ldots,X_1)$. Note that for any $m$ and $c$
there exists $X_m,\ldots ,X_1$ such that
$c=m^{1/\alpha}V_m(X_m,\ldots,X_1)$; the simplest construction is
obtained by letting $X_j=cm^{-1/\alpha}$ for $1 \le j \le m$. Our
suppression of the dependence of $W_n$ on $m$ and $c$ is justified
by Theorem \ref{theorem-Ualpha}, which states that the limiting
value of $W_n$ is the same for all such sequences.

\begin{theorem}
\label{theorem-Ualpha} Let $m \ge 2$ be any integer and suppose
the variables $U_n^{1/\alpha}, \ldots,U_{m+1}^{1/\alpha}, X_m,
\ldots, X_1$ are independent. Let \beq \label{Vnk} V_{n,m} =
V_n(U_n^{1/\alpha}, \ldots,U_{m+1}^{1/\alpha},X_m,\ldots,X_1),
\enq be the optimal two choice value, and suppose
$V_m(X_m,\ldots,X_1)=c \in (0, \infty)$. Then
$$
W_n = n^{1/\alpha}V_{n,m} \quad \mbox{for $n > m$,}
$$
satisfies \beq \label{limWnhca} \lim_{n \rightarrow \infty} W_n =
h(b_\alpha), \enq where $b_\alpha$ is the unique solution to
(\ref{1.13}).

In particular, the optimal two stop value $V_n$ for a sequence of
i.i.d. variables with distribution function ${\cal
U}^\alpha(x)=x^\alpha$ for $0 \le x \le 1$ and $\alpha>0$
satisfies \beq \label{Wnisgood} \lim_{n \rightarrow \infty}n\,
{\cal U}^\alpha(V_n) = h^\alpha(b_\alpha); \enq that is, the
conclusion of Theorem \ref{nancy} holds for the ${\cal U}^\alpha$
family of distributions.
\end{theorem}

\noindent {\bf Proof:} We apply Lemma \ref{nlemma5.2} with the
given $m$ and $c$. Letting $n \rightarrow \infty$ in
(\ref{12star}) and using (\ref{15star}),
$$
d_{j,\alpha} \le \liminf_{n \rightarrow \infty} W_n \le
 \limsup_{n \rightarrow \infty} W_n \le d_\alpha \quad \mbox{for all $j > j_1$.}
$$
Now letting $j \rightarrow \infty$ and using (\ref{n5.6}) gives
(\ref{limWnhca}). The $W_n$ values for the i.i.d. sequence with
distribution function ${\cal U}^{\alpha}$ are generated by
recursion (\ref{21star}) for the particular case $m=2$ and
$c=2^{1/\alpha}E[U_2^{1/\alpha} \wedge U_1^{1/\alpha}]$, thus
proving (\ref{Wnisgood}). $\qed$ \\[1ex]

We conclude this section with some results on the existence of
moments for both the one and two-stop problems.
\begin{theorem}
\label{anmoments} Let $U_n^{1/\alpha}, \ldots,U_1^{1/\alpha}$ be
an i.i.d. sequence with distribution function ${\cal U}^\alpha$,
$a_n$ a sequence of constants in $[0,1]$ with $a_0=1$, and \beq
\label{defTn} T_n = \max \{1 \le k \le n: U^{1/\alpha}_k \le
a_{k-1}\}. \enq When $A_n=n^{1/\alpha} a_n$ satisfies
$$
0< {\underline \kappa}  = \liminf_{n \rightarrow \infty} A_n \le
\limsup_{n \rightarrow \infty} A_n = {\overline \kappa} < \infty
$$
we have \beq \label{finite-m-moment} \limsup_{n \rightarrow
\infty}E(n^{1/\alpha}U^{1/\alpha}_{T_n})^\s < \infty \quad
\mbox{for all $\s<\alpha {\underline \kappa}^\alpha$.} \enq If
$$
{\overline \kappa} < \infty,
$$
we have \beq \label{infinite-m-moment} \lim_{n \rightarrow
\infty}E(n^{1/\alpha}U^{1/\alpha}_{T_n})^\s = \infty \quad
\mbox{for all $\s>\alpha {\overline \kappa}^\alpha$.} \enq

\end{theorem}

\noindent {\bf Proof:} Let
$$
M_n(\s)=E(U^{\s/\alpha}_{T_n})
$$
be the $\s^{th}$ moment of the $a_k$ stopped sequence. The
sequence $M_n(\s)$ satisfies the recursion \beas M_{n+1}(\s) &=&
\int_0^{a_n} x^\s \alpha x^{\alpha -1}dx + (1-a_n^\alpha) M_n(\s),
\quad n \ge 1. \enas Substituting $y=n x^\alpha$, \beas
M_{n+1}(\s) &=& \frac{1}{n^{1+\s/\alpha}} \int_0^{A_n^\alpha}
y^{\s/\alpha} dy + (1-a_n^\alpha) M_n(\s). \enas Multiplying by
$n^{\s/\alpha}$, and letting $n^{\s/\alpha} M_n(\s)=S_n(\s)$,
\bea
\nonumber \left( \frac{n}{n+1} \right)^{\s/\alpha} S_{n+1}(\s) &=&
 \frac{1}{n} \int_0^{A_n^\alpha} y^{\s/\alpha} dy +
       (1-\frac{A_n^\alpha}{n}) S_n(\s) \\
       &=&
       \label{fit}
       \frac{A_n^{\alpha+\s}}{n(1+\s/\alpha)} + (1-\frac{A_n^\alpha}{n})
S_n(\s).
\ena
To show (\ref{finite-m-moment}), multiplying (\ref{fit}) by
$((n+1)/n)^{\s/\alpha}$ and noting that
$$
\left(\frac{n+1}{n} \right)^{\s/\alpha} = 1 + \frac{\s}{\alpha
n}+O_{r/\alpha}(n^{-2})
$$
by the boundedness of the sequence $A_n$ and ${\underline
\kappa}^\alpha > r/\alpha$, we obtain for all $n$ sufficiently
large,
$$
S_{n+1}(\s) \le \frac{2^{\s/\alpha} A_n^{\alpha+\s}}{n(1+\s/\alpha)} +
          (1-\frac{(A_n^\alpha-\s/\alpha)}{2n}) S_n(\s);
$$
(\ref{finite-m-moment}) now follows from Lemma \ref{lemma4.1} A.

To show (\ref{infinite-m-moment}) we note that for all $n$
sufficiently large, using (\ref{fit}),

\beas S_{n+1}(r) &\ge& (\frac{n+1}{n})^{\s/\alpha}
(1-\frac{A_n^\alpha}{n})S_n(r) \\
&=&
 (1 + \frac{\s}{\alpha n} + O_{r/\alpha}(n^{-2}))(1-\frac{A_n^\alpha}{n})S_n(r)\\
 &=&
 \left( 1 + \frac{(\s/\alpha - A_n^\alpha)}{n}+O_{r/\alpha, {\overline \kappa}}(n^{-2})
 \right) S_n(r)\\
&\ge& (1+\frac{(\s/\alpha - A_n^\alpha)}{2n})S_n(r). \enas Now,
recalling that $r/\alpha > {\overline \kappa}^\alpha$, apply Lemma
\ref{lemma4.1} B.

\begin{corollary}
\label{m-moments}
Let $\Optone_n^{U^{1/\alpha}}$ and $\Opttwo_n^{U^{1/\alpha}}$ be the one and two choice
random values obtained from optimally stopping an independent sequence of variables
having distribution ${\cal U}^\alpha$. In the one choice case,
\bea
\label{no-room}
\mbox{if $\s < 1+\alpha$,} \quad
\limsup_{n \rightarrow \infty}E(n^{1/\alpha} \Optone_n^{U^{1/\alpha}})^\s < \infty,\\
\nonumber
\quad \mbox{while if $\s > 1+\alpha$,}
\quad
\limsup_{n \rightarrow \infty}E(n^{1/\alpha} \Optone_n^{U^{1/\alpha}})^\s = \infty.
\ena
In the two choice case,
\beq
\mbox{if $\s < 1+\alpha$,} \quad \limsup_{n \rightarrow \infty}E(n^{1/\alpha}
\Opttwo_n^{U^{1/\alpha}})^\s < \infty.
\enq
\end{corollary}

\noindent {\bf Proof:} For one choice, apply Theorem
\ref{anmoments} with $a_n =V_n^1$, and therefore
$U_{T_n}^{1/\alpha}=\Optone_n^{U^{1/\alpha}}$. By (\ref{1.10}),
$$
\lim_{n \rightarrow \infty} n^{1/\alpha} V_n^1  =
  \lim_{n \rightarrow \infty}  (n \,{\cal U}^\alpha (V_n^1))^{1/\alpha} =
  (1 + 1/\alpha)^{1/\alpha}.
$$
The one choice results now follow from (\ref{finite-m-moment})
and (\ref{infinite-m-moment}) of Theorem \ref{anmoments} with
${\overline \kappa}={\underline \kappa} =
(1+1/\alpha)^{1/\alpha}$.

For two choices, let $T_n$ be defined as in (\ref{defTn}) with
$b_n$, the first choice thresholds given in (\ref{2.3}),
replacing $a_n$, and let $B_n = n^{1/\alpha}b_n$. Then as
$\Opttwo_n^{U^{1/\alpha}} \le U^{1/\alpha}_{T_n}$, it clearly
suffices to show that for $\s < 1+\alpha$,
$$
\limsup_{n \rightarrow \infty} E (n^{1/\alpha}
U^{1/\alpha}_{T_n})^\s < \infty.
$$

Reiterating (\ref{2.9}), $W_n = h_n(B_n^\alpha)$,
and by Theorem \ref{theorem-Ualpha}
$$
\lim_{n \rightarrow \infty} W_n = d_\alpha = h(b_\alpha).
$$

We show $\lim_{n \rightarrow \infty} B_n^\alpha = b_\alpha$.
Suppose $\limsup_{ n\rightarrow \infty} B_n^\alpha = B^\alpha >
b_\alpha$. Then there exists $\epsilon>0$ such that $B^\alpha -
\epsilon > b_\alpha$. But then $\limsup_{n \rightarrow \infty}
h_n(B_n^\alpha) \ge \limsup_{n \rightarrow \infty} h_n (B^\alpha
- \epsilon) =
 h(B^\alpha - \epsilon) > h(b_\alpha)$, a contradiction.
Similarly if $\liminf_{n \rightarrow \infty} B_n^\alpha <
b_\alpha$. Thus the limit of $B_n$ exists and
$$
n^{1/\alpha}b_n = B_n \rightarrow b_\alpha^{1/\alpha}.
$$
By (\ref{finite-m-moment}) it suffices to show that $b_\alpha > 1
+ 1/\alpha$, which, by Lemmas \ref{nlemma4.2} and
\ref{nlemma5.1}, would follow from $H(1+1/\alpha)>0$. Now \beas
H(1+1/\alpha) &=& (1+\alpha)^{1/\alpha}
    \left[ \int_0^{1+1/\alpha} \left( \frac{y}{\alpha + 1 + \alpha y} \right)^{1/\alpha}dy
       - \left( \frac{1+1/\alpha}{2(1+\alpha)}  \right)^{1/\alpha}  \right] \\
&>& (1+\alpha)^{1/\alpha}
    \left[ \left(\frac{1}{2(1+\alpha)}\right)^{1/\alpha}\int_0^{1+1/\alpha} y^{1/\alpha}dy
       - \left( \frac{1+1/\alpha}{2(1+\alpha)}  \right)^{1/\alpha}  \right]\\
&=& 0,
\enas
completing the proof. $\qed$

\begin{remark} Kennedy and Kertz (1991, Theorem 1.4) obtain the limiting distribution of
the scaled optimal one stop random variable
$n^{1/\alpha}{\bf 1}_n^{U^{1/\alpha}}$. It is easily checked that this
limiting distribution has a finite $\s^{th}$ moment if and only if $\s < 1 +
\alpha$, which is not surprising, when compared with (\ref{no-room}) in
Corollary \ref{m-moments}.
\end{remark}

\begin{remark}
\From the proof that $b_\alpha > 1 +1/\alpha$ in Corollary
\ref{m-moments}, it follows that the limiting thresholds $b_n$ for
the first choice in the optimal two-choice problem are larger than
the corresponding values $V_n^1$ for the optimal one choice
problem, for all $\alpha > 0$. This is reasonable, as with two
choices one `can afford' to make the first of the two choices in
the two stop problem earlier than the only choice in the one stop
problem.

Another interpretation of the inequality $b_n > V_n^1$ is gained
by applying $V_n^1( )$ to both sides, to obtain $V_n > V_{2n}^1$
i.e. one is better off having one choice among $2n$ variables
than having two choices among $n$ variables.
\end{remark}

\begin{remark} Whereas it follows from Resnick, (1987,
Proposition 2.1) that all scaled moments of the minimum exist, it is of
interest to note that no moment with $\s > 1 + \alpha$  exists for the
optimal scaled one-choice value.
\end{remark}

\noindent
\section{Extension to General Distributions}
Theorem \ref{theorem-Ualpha} treats the special case where the
variables have distribution function ${\cal U}^{\,\alpha}(x)$ as
in (\ref{1.14}). At the end of this section we prove Theorem
\ref{nancy} for an i.i.d. sequence of random variables in a much
wider class.

To prove Theorem \ref{nancy} the two stop problem is considered
for $X_n, \ldots, X_1$, non-trivial independent but not
necessarily identically distributed random variables. It is direct
to see that the dynamic programming equations given in the
introduction for an i.i.d. sequence hold under the assumption of
independence alone. In particular, the one and two stop value
functions $V_n^1(x)$ and $V_n^2(x)$ are again given through
(\ref{1.7}) and (\ref{1.9}) respectively. With nothing guaranteed,
we have that $V_n^1=V_n^1(\infty)$ and $V_n=V_n^2(\infty)$ are the
one and two choice optimal stopping values, respectively. However,
Lemma \ref{bkexists} gives an alternative representation for
$V_n^1$ which reduces to $V_n^1(x_F)$ as given earlier for the
i.i.d. case, as well as conditions which guarantee that the
`threshold' indifference sequences are uniquely defined for
independent but not necessarily identically distributed sequences.

\begin{lemma} \label{bkexists} Let $X_n,
\ldots, X_1$ be non-negative independent random variables with
distribution functions $F_n,\ldots,F_1$ respectively, and $x_F$
given in (\ref{defxF}). Then for all $x \ge 0$ the function
$V_k^1(x)$ given by (\ref{1.7}) satisfies $0 \le V_k^1(x) \le x$
and is non-decreasing and continuous.

Letting $v_1=x_{F_1}$ and
\beas
v_k=v_{k-1} \wedge w_k \quad \mbox{where $w_k = \inf\{y:
V_{k-1}^1(y) \ge x_{F_k}\}$} \quad \mbox{for $2 \le k \le n$,}
\enas
the function $V_k^1(x)$ is strictly monotone increasing for $0 < x
< v_k$, and satisfies $V_k^1(x)=V_k^1(v_k)$ for $x \ge v_k$; in
particular
$$
V_k^1=V_k^1(v_k).
$$

Furthermore, the indifference numbers $b_k,\, 2 \le k \le n$ given
by the solutions to
$$
V_k=V_k^1(b_k)
$$
exist, and are uniquely defined in $[0,v_k]$.
\end{lemma}
{\bf Proof:} For all $x \ge 0$ the function $V_1^1(x)=E[X_1 \wedge
x]=x-\int_0^xF_1(y)dy$ satisfies $0 \le V_1^1(x) \le x$ and is
non-decreasing and continuous; further, $V_1^1(x)$ is strictly
increasing for $x$ in $[0,x_{F_1}]$, and $V_1^1(x)=EX_1$ for $x
\ge x_{F_1}$. Now assume for all $x \ge 0$ that $0\le
V_{k-1}^1(x)\le x$, and $V_{k-1}^1(x)$ is non-decreasing and
continuous. Then $V_k^1(x) = E[V_{k-1}^1(x) \wedge X_k] \le E[x
\wedge X_k] \le x$. For $0 \le x \le y$ we have $V_k^1(x) =
E[V_{k-1}^1(x) \wedge X_k] \le E[V_{k-1}^1(y) \wedge X_k] =
V_k^1(y)$, and $V_k^1(x)$ is continuous for all $x$ by the bounded
convergence theorem, using the continuity of $V_{k-1}^1(x)$ and
its upper bound of $x$.

To prove strict monotonicity, assume that $V_{k-1}^1(x)$ is
strictly increasing for $0 \le x \le v_{k-1}$ and take $0 \le x <
y \le v_k$. Since $v_k \le v_{k-1}$ we have
$V_{k-1}^1(x)<V_{k-1}^1(y)$, and since $x<w_k$ we have
$V_{k-1}^1(x)<x_{F_k}$ and therefore $P(X_k>V_{k-1}^1(x))>0$.
Hence
$$
V_k^1(x) = E[V_{k-1}^1(x) \wedge X_k] < E[V_{k-1}^1(y) \wedge X_k] = V_k^1(y).
$$
Now, assuming that $V_{k-1}^1(x)$ is constant for $x  \ge
v_{k-1}$, then for all $x \ge v_{k-1}$,
\beas
V_k^1(x) &=& E[V_{k-1}^1(x) \wedge X_k] = E[V_{k-1}^1(v_{k-1})
\wedge X_k] =
E[V_{k-1}^1(v_{k-1}) \wedge x_{F_k} \wedge X_k]\\
&=& E[V_{k-1}^1(v_{k-1}) \wedge V_{k-1}^1(w_k) \wedge
X_k]=E[V_{k-1}^1(v_{k-1} \wedge w_k) \wedge X_k]= V_k^1(v_k).
\enas
Similarly, for all $x \ge w_k$,
\beas
V_k^1(x) = E[V_{k-1}^1(x) \wedge X_k]=E[V_{k-1}^1(x) \wedge
V_{k-1}^1(w_k) \wedge X_k] = E[V_{k-1}^1(w_k) \wedge
X_k]=V_k^1(w_k),
\enas
from which it follows that $V_k^1(v_k)=V_k^1(w_k)$ and
$V_k^1(x)=V_k^1(v_k)$ for all $x \ge v_k$.

Since
$$
0 = V_k^1(0) \le V_k=  V_k^1(b_k)  \le V_k^1 = V_k^1(v_k),
$$
and $V_k^1(x)$ is continuous and strictly monotone increasing in
$[0,v_k]$, the solution $b_k$ exists uniquely in $[0,v_k]$.$\qed$

In the case where the variables are i.i.d., since $V_{k-1}^1(y)
\le y$ we have $w_k \ge x_F$, and hence $v_k=x_F$, as given in
Section \ref{intro}.

\begin{lemma} \label{bnmonotone} For any sequence of nonnegative independent
random variables $X_n,\ldots,X_1$ the sequence $b_k,2 \le k \le n$
is monotone non-increasing.
\end{lemma}

\noindent {\bf Proof:} We first show that
$$
V_{k+1}^2 \le E[X_{k+1} \wedge V_k^2], \quad \mbox{$2 \le k \le
n-1$.}
$$
The right hand side is the value obtained by applying, on the
sequence $X_{k+1},\ldots,X_1$, the suboptimal two choice rule
where $X_{k+1}$ is chosen as the first and second choice if
$X_{k+1}<V_k^2$ (this is the same as taking $X_{k+1}$ as the first
choice and not taking any second choice), and when $X_{k+1}\ge
V_k^2$ the optimal two choice rule is applied on $X_k,\ldots,X_1$.
The inequality reflects that the optimal rule does as well as
this, or any other, two choice rule on this sequence. Therefore
$$
V_{k+1}^1(b_{k+1}) = V_{k+1}^2 \le E[X_{k+1} \wedge V_k^2] =
E[X_{k+1} \wedge V_k^1(b_k)] = V_{k+1}^1(b_k).
$$
Since $V_{k+1}^1(x)$ is strictly monotone increasing in the
interval $[0,v_{k+1}]$, which contains $b_{k+1}$, the Lemma is
shown. $\qed$

Below we consider stochastic dominance between two random
variables, and write $Y \le_d X$ when $P(Y>t) \le P(X>t)$ for all
$t$.
\begin{lemma}
\label{taulemma} Let $X_n,\ldots,X_1$ and $Y_n,\ldots,Y_1$ be
sequences of independent non-negative random variables having two
choice value and threshold sequences $V^X_j,V^Y_j$ and
$b_j^X,b_j^Y, j=1,\ldots, n$ respectively. If for some $m \ge 2$,
\bea \label{j=12} Y_j \le_d X_j, \quad j=1,\ldots,m, \ena and there
exists $\tau\ge \max\{b_m^X, b_m^Y\}$ such that
\bea \label{jlarge}
\tau \wedge Y_{j+1} \le_d \tau \wedge X_{j+1} \quad
 \mbox{for $m \le j < n$},
\ena
then
\bea
\label{VjYleVjX}
V_j^Y \le V_j^X, \quad \mbox{for $j=2,\ldots,n$;}
\ena
hence, if the inequalities in (\ref{j=12}) and (\ref{jlarge}) are
replaced by equalities, then $V_j^Y = V_j^X, j=2,3,\ldots,n$.
Finally, $V_n^X$ is unchanged upon replacing any $X_{j+1}$ by
$\tau \wedge X_{j+1}, 2 \le j < n$, for any $\tau \ge b_j^X$.
\end{lemma}

\noindent {\bf Proof:} Let $V_n^{X,1}(x)$ and $V_n^{Y,1}(x)$
denote the optimal one choice value functions for the $X$ and $Y$
sequences respectively, with guaranteed value $x$, as in
(\ref{1.7}). A simple induction using (\ref{j=12}) gives
$V_j^{Y,1}(x) \le V_j^{X,1}(x)$ for all $x$ and $1 \le j \le m$.

First suppose that (\ref{jlarge}) holds for some arbitrary $\tau$,
and that for some $m \le j < n$,
\bea
\label{VjYdominated}
V_j^{Y,1}(x) \le V_j^{X,1}(x) \quad \mbox{for
all $x \le \tau$.}
\ena
Then for $x \le \tau$, using that $V_j^{X,1}(x) \le x \le \tau$
and $V_j^{Y,1}(x) \le x \le \tau$ by Lemma \ref{bkexists}, we have
$$
Y_{j+1} \wedge V_j^{Y,1}(x) \le Y_{j+1} \wedge V_j^{X,1}(x) =
(Y_{j+1} \wedge \tau) \wedge V_j^{X,1}(x) \le_d (X_{j+1} \wedge
\tau ) \wedge V_j^{X,1}(x) = X_{j+1} \wedge V_j^{X,1}(x),
$$
giving \beas V_{j+1}^{Y,1}(x) = E[Y_{j+1} \wedge V_j^{Y,1}(x)]
 \le E[X_{j+1} \wedge V_j^{X,1}(x)]
 = V_{j+1}^{X,1}(x), \quad \mbox{for $x \le \tau,$}
 \enas
and thus (\ref{VjYdominated}) holds for $1 \le j \le n$.

For $\tau \ge \max\{b_j^X,b_j^Y\}$ and $j=m$, Lemma
\ref{bnmonotone} implies this inequality holds for $m \le j <n$,
and therefore, for instance,
$$
Y_{j+1} \wedge b_j^Y = (Y_{j+1} \wedge \tau) \wedge b_j^Y \le_d
(X_{j+1} \wedge \tau) \wedge b_j^Y = X_{j+1} \wedge b_j^Y.
$$
Now note that (\ref{j=12}) yields $V_j^Y \le V_j^X$ for $1 \le j
\le m$, so assuming this inequality for some $j$, $m \le j <n$, we
now have
\beas V_{j+1}^Y
&=& E[V_j^{Y,1}(Y_{j+1}) \wedge V_j^Y] = E[V_j^{Y,1}(Y_{j+1})
\wedge V_j^{Y,1}(b_j^Y) \wedge V_j^Y]=
E[V_j^{Y,1}(Y_{j+1} \wedge b_j^Y) \wedge V_j^Y]\\
&\le& E[V_j^{Y,1}(X_{j+1} \wedge b^Y_j) \wedge V_j^Y]
\le E[V_j^{X,1}(X_{j+1} \wedge b^Y_j) \wedge V_j^Y] \\
&=& E[V_j^{X,1}(X_{j+1}) \wedge V_j^{X,1}(b_j^Y) \wedge
V_j^{Y,1}(b^Y_j) ]=
E[V_j^{X,1}(X_{j+1}) \wedge V_j^{Y,1}(b^Y_j)]\\
&=& E[V_j^{X,1}(X_{j+1}) \wedge V_j^Y] \le E[V_j^{X,1}(X_{j+1}) \wedge V_j^X]\\
&=& V_{j+1}^X; \enas thus (\ref{VjYleVjX}) holds.

To prove the final claim, let $Y_n,\ldots,Y_1$ be the sequence
where any number of variables $X_{j+1},  2 \le j < n$ have been
replaced by $X_{j+1} \wedge \tau$ with $\tau \ge b_j^X$. Note that
(\ref{j=12}) and (\ref{jlarge}) hold with equality, and hence so
does (\ref{VjYdominated}). Clearly $V_2^Y=V_2^X$, so assuming
$V_j^Y=V_j^X$ for $2 \le j <n$, we have, taking the non-trivial
case of $j$ for which $Y_{j+1}=X_{j+1} \wedge \tau$
\beas
V_{j+1}^Y&=&E[V_j^{Y,1}(Y_{j+1}) \wedge V_j^Y]=E[V_j^{Y,1}(X_{j+1}
\wedge \tau) \wedge V_j^Y]=E[V_j^{X,1}(X_{j+1}) \wedge
V_j^{X,1}(\tau) \wedge V_j^X]\\
&=& E[V_j^{X,1}(X_{j+1}) \wedge V_j^{X,1}(\tau) \wedge
V_j^{X,1}(b_j^X)]=E[V_j^{X,1}(X_{j+1}) \wedge
V_j^{X,1}(b_j^X)]=V_{j+1}^X. \qed
\enas

Let now $X_n,\ldots,X_1$ be i.i.d. as $X$ with distribution
function $F$ satisfying the hypotheses of Theorem \ref{nancy}.
Without loss of generality we may assume that the function $L$ in
(\ref{wlog1}) satisfies $\lim_{x \downarrow 0}L(x)=1$, since if
$F_X(x)=x^\alpha L_\B(x)$ with $\lim_{x \downarrow 0}L_\B(x)=\B
\in (0,\infty)$, then $Z=\B^{1/\alpha}X$ has distribution function
$F_Z(z)=z^\alpha (1/\B) L_\B(\B^{-1/\alpha}z)$ with $\lim_{z
\downarrow 0}(1/\B)L_\B(\B^{-1/\alpha}z)=1$. Since $V_n^Z =
\B^{1/\alpha} V_n^X$, we have
$$
F_Z(V_n^Z) = F_X(V_n^X),
$$
and hence we can assume that $X$ has distribution function $F$
such that \beq \label{dontcare} F(x) = x^\alpha L(x) \quad
\mbox{$\lim_{x \downarrow 0}L(x) = 1.$} \enq

\begin{corollary}
\label{Xwlogbdd} Let $X_n,\ldots,X_1$ be a sequence of i.i.d.
non-negative random variables with $E[X_2 \wedge X_1]<\infty$ and
distribution function satisfying (\ref{dontcare}). Then there
exists an i.i.d. sequence $Y_n,\ldots,Y_1$ of bounded non-negative
random variables with distribution function satisfying
(\ref{dontcare}) such that $V_n^Y=V_n^X$ for all $n \ge 2$.
\end{corollary}

\noindent {\bf Proof:} Assume $x_F=\infty$, else there is nothing
to prove. For all $x >0$ sufficiently small, using the
non-degeneracy of the distribution $F$ on $[0,x]$, Jensen's
inequality applied to the concave function $\psi(u)= u \wedge x$
yields
$$
E[x \wedge X_1] \le x \wedge EX_1, \quad \mbox{with strict
inequality for all $x$ sufficiently small.}
$$
Thus
$$
E[X_2 \wedge X_1|X_2] \le X_2 \wedge EX_1, \quad \mbox{with strict inequality having positive probability}
$$
and therefore, since $V_1^1(\infty) = EX_1$ (which may be infinite),
$$
0< V_2= E[X_2 \wedge X_1] < E[X_2 \wedge EX_1] = V_2^1(\infty).
$$
Using $x_F=\infty$ and Lemma \ref{bkexists}, $V_2^1(x)$ is
continuous and strictly monotone increasing on $(0,\infty)$, hence
the solution $b_2$ to
$$
V_2=V_2^1(x)
$$
exists, is unique, and satisfies $0<b_2 < \infty$.

For $j=1,\ldots,n$ and any $K \ge b_2$ let \beas Y_j = \left\{
\begin{array}{cl}
X_j & \mbox{for $X_j \le b_2$}\\
K  & \mbox{for $X_j > b_2.$}
\end{array} \right.
\enas
Using Lemma \ref{taulemma} with $m=2$ and $\tau = b_2$, we see
that the two stop values of $X_j,\ldots,X_1$ and of
$Y_j,\ldots,Y_3,X_2,X_1$ are the same for $2 \le j \le n$, i.e.
$V_j^X = V_j^2(X_j,\ldots,X_1) = V_j^2(Y_j,\ldots,Y_3,X_2,X_1)$.
Since the distribution of $X_j$ is unbounded, $P(X_j
> b_2)>0$, which guarantees that $K \ge b_2$ can be chosen to
yield $E[Y_2 \wedge Y_1] = E[X_2 \wedge X_1]$. But now, with the
equality $V_j^2(Y_j,\ldots,Y_1)=V_j^2(Y_j,\ldots,Y_3,X_2,X_1)$
between the optimal two stop values on the sequences indicated now
true for $j = 2$ by choice of $K$, assuming it true for $j\ge 2$
and using the notation as in the proof of Lemma \ref{taulemma}
yields
\beas
V_{j+1}^Y&=&V_{j+1}^2(Y_{j+1},\ldots,Y_1)=E[V_j^{Y,1}(Y_{j+1})
\wedge V_j^2(Y_j,\ldots,Y_1)]\\
&=&E[V_j^{Y,1}(Y_{j+1}) \wedge
V_j^2(Y_j,\ldots,Y_3,X_2,X_1)] =
V_{j+1}^2(Y_{j+1},\ldots,Y_3,X_2,X_1)=V_{j+1}^X.
\enas
Since $P(X_j \le x)=P(Y_j \le x)$ for all $0 \le x < b_2$, the
distribution $P(Y_j \le x)$ satisfies (\ref{dontcare}) and the
bounded i.i.d. sequence $Y_n,\ldots,Y_1$
has all the claimed properties.$\qed$ \\[1ex]

We have the following Lemma.
\begin{lemma}
\label{Linverse} Let $X$ have distribution function $F(x)=P(X \le
x)$, and set \beas F^{-1}(u) = \sup \{ x: F(x) < u \} \quad
\mbox{for $0<u<1$}. \enas Then \beq \label{sets=} F(x) \ge u \quad
\mbox{if and only if} \quad x \ge F^{-1}(u), \enq and with $U
\sim {\cal U}(0,1)$ we have \beq \label{pit} X =_d F^{-1}(U). \enq

In addition, if
$$
F(x)=x^\alpha L_F(x), \quad \mbox{for all $x \ge 0$,
with $\lim_{x \downarrow 0}L_F(x) = 1$,}
$$
then there exists a function $L^*$ such that
\beq
\label{Finverseslow}
F^{-1}(u) = u^{1/\alpha}L_{F^{-1}}(u) =
 u^{1/\alpha}L^*(u^{1/\alpha}), \quad \mbox{with} \quad
\lim_{u \downarrow 0}L^*(u)=1, \enq so that by (\ref{pit}) and
(\ref{Finverseslow}), \beq \label{94} X =_d
U^{1/\alpha}L^*(U^{1/\alpha}). \enq

\end{lemma}
\noindent {\bf Proof:} Let $A_u =\{x:F(x)<u\}$. If $F(x) \ge u$
then $x \not \in A_u$ and therefore $F^{-1}(u) \le x$. If
$F(x)<u$ then by right continuity there exists $\epsilon>0$ such
that $F(x+\epsilon)<u$. Thus $x+\epsilon \in A_u$, which gives
that $F^{-1}(u) \ge x +\epsilon > x$. This demonstrates
(\ref{sets=}). Now replacing $u$ by a random variable $U$ having
the ${\cal U}[0,1]$ distribution we obtain (\ref{pit}), by
$P(F^{-1}(U)\le x)=P(U \le F(x))=F(x)$.

The claim in (\ref{Finverseslow}) is equivalent to \beq
\label{lim=1} \lim_{u \downarrow 0}\frac{(F^{-1}(u))^\alpha}{u} =
1. \enq Using that $F(x)=x^\alpha L_F(x)$,
$$
F^{-1}(u) = \sup \{ x: x^\alpha L_F(x) < u \},
$$
and hence, setting $L_\alpha(y) = L_F(y^{1/\alpha})$,
\beq
\label{Finv-forms}
(F^{-1}(u))^\alpha = \sup \{ x^\alpha: x^\alpha L_F(x) < u \}
                   = \sup \{ y: y L_F(y^{1/\alpha}) < u \}
                   = \sup \{ y: y L_\alpha(y) < u \}.
\enq
Note $y L_\alpha(y) = F(y^{1/\alpha})$ is non-decreasing. Let $\epsilon \in (0,1)$ be
given. Since $\lim_{y \downarrow 0}L_\alpha(y)=1$, there exists $\delta>0$ such that
\beq
\label{Lbounds}
1-\epsilon < L_\alpha(y) < 1+ \epsilon \quad \mbox{for all $0< y < \delta$.}
\enq
Let $0<u<\delta (1-\epsilon)$. Then if $0<y<u/(1+\epsilon)$ we have $y < \delta$ and so
$$
yL_\alpha(y) < y (1+\epsilon) < u,
$$
so
$$
\{y: 0<y<u/(1+\epsilon) \} \subset \{y:yL_\alpha(y)<u \}.
$$
Thus
$$
u/(1+\epsilon) \le (F^{-1}(u))^\alpha \quad \mbox{for all $0< u <
\delta(1-\epsilon)$}.
$$
Now, with $0<u<\delta (1-\epsilon)$ and any $y \in (u/(1 -
\epsilon), \delta)$, by (\ref{Lbounds}),
\beas
%\label{ytoobig}
u<(1-\epsilon)y < yL_\alpha(y),
\enas
and it follows by (\ref{Finv-forms}) that
$$
(F^{-1}(u))^\alpha \le u/(1-\epsilon).
$$
Hence,
$$
1/(1+\epsilon) \le \frac{(F^{-1}(u))^\alpha}{u} \le 1/(1-\epsilon)
\quad \mbox{for $0<u<\delta (1-\epsilon)$,}
$$
and (\ref{lim=1}) is shown. $\qed$

\begin{lemma}
\label{ui} Let $\X_n, n=1,2,\ldots $ be a uniformly integrable
non-negative sequence of random variables, and $0 \le L_n \le L$,
$L$ a constant, with $L_n \rightarrow_p 1$ as $n \rightarrow
\infty$. Then \beas \limsup_{n \rightarrow \infty} E\X_n L_n &=&
\limsup_{n \rightarrow \infty} E\X_n \enas so that in particular,
if $\lim_{n \rightarrow \infty} E\X_n$ exists,
$$
\limsup_{n \rightarrow \infty} E\X_nL_n = \lim_{n \rightarrow \infty} E\X_n.
$$
\end{lemma}

\noindent {\bf Proof:} Let $\epsilon>0$ be given. Since $\X_n$ is uniformly integrable,
there exists $\delta >0$ such that
\beq
\label{Xnu.i.}
E\X_n{\bf 1}_A \le \epsilon \quad \mbox{whenever} \quad P(A) \le \delta.
\enq
Since $L_n \rightarrow_p 1$ as $n \rightarrow \infty$, there exists $n_0$ such that
for all $n \ge n_0$
$$
\Omega_n = \{ |L_n-1| \le \epsilon \} \quad \mbox{satisfies} \quad
P(\Omega_n) \ge 1 - \delta.
$$
Hence, for $n \ge n_0$, using (\ref{Xnu.i.}) and that $\X_n \ge 0$,
with $A= \Omega_n^c$,
\beas
(1-\epsilon) E\X_n {\bf 1}_{\Omega_n} \le E\X_n L_n
  \le (1+\epsilon)E\X_n {\bf 1}_{\Omega_n} + L \epsilon \le (1+\epsilon)E\X_n + L \epsilon
\enas
and
\beas
E\X_n - \epsilon \le E\X_n{\bf 1}_{\Omega_n},
\enas
so that for $n \ge n_0$ we have
$$
(1-\epsilon) (E\X_n - \epsilon) \le E\X_n L_n \le (1+\epsilon)E\X_n +L \epsilon.
$$
Taking $\limsup$ and recalling $\epsilon>0$ was arbitrary
completes the proof. $\qed$ \\[2ex]

\begin{lemma}
\label{usareui}
Let $X_n,\ldots,X_1$ be an integrable i.i.d. sequence with
distribution function
$F(x)$ satisfying (\ref{dontcare}). Let $W_n^X=n^{1/\alpha}V_n^X$
and $W_n^{U^{1/\alpha}}=n^{1/\alpha}V_n^{U^{1/\alpha}}$. Then
$$
\limsup_{n \rightarrow \infty} W_n^X \le \lim_{n \rightarrow \infty} W_n^{U^{1/\alpha}}.
$$
\end{lemma}
\noindent {\bf Proof:} Using Lemma \ref{Linverse}, we construct
i.i.d. pairs $(U_i,X_i)$ with $U_i \sim {\cal U}, X_i \sim F$, and
$$
X_i=U_i^{1/\alpha}L^*(U_i^{1/\alpha}).
$$

By Corollary \ref{Xwlogbdd}, without loss of generality we
can take the $X$ variables to be bounded, and since $L^*(u) \rightarrow 1$ as
$u \downarrow 0$, it follows that $L^*$ is bounded.

Let $\Opttwo_n^{U^{1/\alpha}}$ and  $\Opttwo_n^X$ be the optimal
random $n$-variable two-stop value for the
$U_n^{1/\alpha},\ldots,U_1^{1/\alpha}$ and $X_n,\ldots,X_1$
sequences respectively. Since $E
n^{1/\alpha}\Opttwo_n^{U^{1/\alpha}}=n^{1/\alpha}V_n^{U^{1/\alpha}}
= W_n^{U^{1/\alpha}}$ converges (to $h(b_\alpha)$), we have
$$
P(\Opttwo_n^{U^{1/\alpha}} > \epsilon) =
P(n^{1/\alpha}\Opttwo_n^{U^{1/\alpha}} > n^{1/\alpha} \epsilon) \le
\frac{W_n^{U^{1/\alpha}}}{n^{1/\alpha}\epsilon} \rightarrow 0 \quad
\mbox{as $n \rightarrow \infty$.}
$$
Hence $\Opttwo_n^{U^{1/\alpha}} \rightarrow_p 0$, and therefore
$L^*(\Opttwo_n^{U^{1/\alpha}}) \rightarrow_p 1$. Furthermore,
by Corollary \ref{m-moments}, the collection
$n^{1/\alpha}\Opttwo_n^{U^{1/\alpha}}$ has a bounded $\s^{th}$ moment for some $\s>1$
and hence is uniformly integrable.

Let $\Opttwo_n^{X,U^{1/\alpha}}$
denote the $X$ sequence stopped on the optimal rules
for the $U^{1/\alpha}$ sequence.
Then $\Opttwo_n^{X,U^{1/\alpha}}=
\Opttwo_n^{U^{1/\alpha}}L^*(\Opttwo_n^{U^{1/\alpha}})$, and since these
rules may not be optimal for the $X$ sequence we have
$$
E n^{1/\alpha}\Opttwo_n^X \le E n^{1/\alpha}\Opttwo_n^{X,U^{1/\alpha}}
= E n^{1/\alpha} \Opttwo_n^{U^{1/\alpha}}L^*(\Opttwo_n^{U^{1/\alpha}}).
$$
Taking limsup and using that
$n^{1/\alpha}\Opttwo_n^{U^{1/\alpha}}$ is uniformly integrable
and $L^*$ is bounded and $L^*(\Opttwo_n^{U^{1/\alpha}})
\rightarrow_p 1$, the result follows from Lemma \ref{ui} and the
fact that $W_n^{U^{1/\alpha}}$ converges. $\qed$

\begin{lemma}
\label{bn->0}
 Let $X_n,\ldots,X_1$ be i.i.d. random variables with distribution
function $F$ satisfying (\ref{dontcare}). Then the indifference values $b_n$
for $X$ satisfy
$$
\lim_{n \rightarrow \infty} b_n =0.
$$
\end{lemma}
\noindent {\bf Proof:} Let $V_n^1(x)=V_n^1(X_n,\ldots,X_1;x)$ and
$V_n^1(X_n,\ldots,X_1)$ denote the optimal one stop value on
$X_n,\ldots,X_1$ with and without the guaranteed bound of $x$,
respectively. Note that trivially for $k=1$ we have that
$$
V_k^1(X_k,\ldots,X_1;x)=V_k^1(X_k \wedge x,\ldots,X_1 \wedge x),
$$
and assuming it true for some $k$, $1\le k<n$ and using
$V_k^1(X_k,\ldots,X_1;x) =V_k^1(x) \le x$ gives
\beas
&& V_{k+1}^1(X_{k+1},\ldots,X_1;x)=E[X_{k+1} \wedge
V_k^1(X_k,\ldots,X_1;x)]=E[(X_{k+1} \wedge x) \wedge
V_k^1(X_k,\ldots,X_1;x)]\\
&=& E[(X_{k+1} \wedge x) \wedge V_k^1(X_k \wedge x,\ldots,X_1
\wedge x)] =V_{k+1}^1(X_{k+1} \wedge x,\ldots,X_1 \wedge x).
\enas

Since $b_n$ is monotone non-increasing by Lemma \ref{bnmonotone},
$b_n \downarrow b \ge 0$, and we have \beq \label{doescare}
V_n^1(X_n \wedge b, \ldots, X_1 \wedge b) = V_n^1(b) \le
V_n^1(b_n) = V_n^X. \enq Hence the two choice value $V_n^X$ on
$X_n,\ldots,X_1$ is greater (worse) than the optimal one choice
value of the sequence of i.i.d. random variables $b \wedge
X_n,\ldots,b \wedge X_1$. If $b>0$, by (\ref{1.10}), the limit of
the scaled optimal one choice value of this sequence, $W_n^{X
\wedge b,1}$ say, is the same as the limit of $W_n^{X,1}$, the
scaled optimal one choice value for $X_n,\ldots,X_1$. But then,
using (\ref{doescare}) in the first inequality, Lemma
\ref{usareui} for the second inequality, Theorem
\ref{theorem-Ualpha} for the equality, (\ref{2<1}) for the strict
inequality and the results of Kennedy and Kertz (1991) for the
last two equalities we have
$$
\lim_{n \rightarrow \infty } W_n^{X,1} \le
 \limsup_{n \rightarrow \infty } W_n^X \le
 \lim_{n \rightarrow \infty}W_n^{U^{1/\alpha}} = h(b_\alpha) < (1+1/\alpha)^{1/\alpha}
 =\lim_{n \rightarrow \infty}W_n^{U^{1/\alpha},1}
 =\lim_{n \rightarrow \infty} W_n^{X,1},
$$
a contradiction. $\qed$

\begin{lemma}
\label{Vnn} Let $(U_i,X_i),\,\, i=n,\ldots,1$ be independent pairs
of random variables with $U_i$ uniform on $[0,1]$ and $X_i$ having
distribution function $F$ satisfying (\ref{dontcare}). Let
$V_{n,m}$ be defined as in (\ref{Vnk}), giving in particular
$V_{n,n}=V_n^X$. Then for every $\epsilon \in (0,1)$, there exists
$m$ such that \bea \label{Vnkinfsup} \frac{1}{1+\epsilon}  \le
\liminf_{n \rightarrow \infty} \frac{V_{n,m}}{V_{n,n}} \le
\limsup_{n \rightarrow \infty}\frac{V_{n,m}}{V_{n,n}}
  \le \frac{1}{1-\epsilon}.
\ena
\end{lemma}

\noindent {\bf Proof:} Using (\ref{94}) of Lemma \ref{Linverse},
we can construct the i.i.d. $X$ sequence using an i.i.d. sequence
$U^{1/\alpha}$ with distribution ${\cal U}^\alpha$ by defining
$X_i$ as
\beq
\label{YgeXL(X)2} X_i = U^{1/\alpha}_i L^*(U_i^{1/\alpha}) \quad
\mbox{a.s.} \enq where $\lim_{u \downarrow 0}L^*(u)=1$.
Hence, for the given $\epsilon \in (0,1)$ there exists $\delta>0$ such that
\beq
\label{LFXbound}
1-\epsilon \le L^*(u^{1/\alpha}) \le 1+\epsilon \quad \mbox{for
$0<u \le \delta$,}
\enq
and so by (\ref{YgeXL(X)2}) and (\ref{LFXbound}) we have
$$
(1+\epsilon)^{-1}X_i  \le
 U_i^{1/\alpha} \le
 (1-\epsilon)^{-1}X_i \quad \mbox{when $U_i \le \delta$}.
$$

By condition (\ref{dontcare}), $F$ is continuous at 0 and satisfies $F(0)=0$,
and therefore there exists $\rho>0$ with $0< F(\rho) \le \delta$. But by
(\ref{sets=}), since
$$
U_i \le F(\rho) \quad \mbox{if and only if} \quad  X_i \le \rho,
$$
we have
$$
\mbox{if} \quad X_i \le \rho \quad \mbox{then} \quad U_i \le \delta.
$$

Let $\tau = \min\{ \delta, \rho\}$, and $b_n^X$ and
$b_n^{U^{1/\alpha}}$ be the indifference values for the $X$ and
$U^{1/\alpha}$ variables, respectively, which by Lemma
\ref{bn->0} converge monotonically to zero. Hence there exists
$m$ with
$$
\max\{ b_m^{U^{1/\alpha}}, b_m^X\} \le \tau,
$$
and for all $n \ge m$, by Lemma \ref{taulemma}, \beas
&&(1+\epsilon)^{-1} V_n(X_n,\ldots,X_1)\\
  &=& (1+\epsilon)^{-1} V_n(X_n \wedge \tau,\ldots,X_{m+1} \wedge \tau, X_m, \ldots, X_1)\\
&=& V_n((1+\epsilon)^{-1} (X_n \wedge \tau),\ldots,
    (1+\epsilon)^{-1} (X_{m+1} \wedge \tau),
 (1+\epsilon)^{-1}X_m ,\ldots,(1+\epsilon)^{-1} X_1)\\
&\le& V_n(U_n^{1/\alpha} \wedge \tau,\ldots,
    U_{m+1}^{1/\alpha} \wedge \tau, X_m,\ldots,X_1)\\
&=& V_n(U_n^{1/\alpha},\ldots, U_{m+1}^{1/\alpha}, X_m,\ldots,X_1)\\
&=& V_n(U_n^{1/\alpha} \wedge \tau,\ldots,
    U_{m+1}^{1/\alpha} \wedge \tau, X_m,\ldots,X_1)\\
&\le& V_n((1-\epsilon)^{-1}(X_n \wedge \tau),\ldots,
    (1-\epsilon)^{-1}(X_{m+1} \wedge \tau),
  (1-\epsilon)^{-1}X_m,\ldots,(1-\epsilon)^{-1}X_1)\\
&\le& (1-\epsilon)^{-1} V_n(X_n \wedge \tau,\ldots,
    X_{m+1} \wedge \tau, X_m,\ldots,X_1)\\
&=& (1-\epsilon)^{-1} V_n(X_n,\ldots,X_1). \enas Now dividing by
$V_{n,n}$ we see that for all $n \ge m$,
$$
\frac{1}{1+\epsilon} \le \frac{V_{n,m}}{V_{n,n}} \le
\frac{1}{1-\epsilon},
$$
completing the proof. $\qed$ \\[1ex]

\noindent {\bf Proof of Theorem \ref{nancy}:} Clearly, for all $0
\le m \le n$,
\beas
\frac{V_n(U_n^{1/\alpha},\ldots,U_1^{1/\alpha})}{V_n(X_n,\ldots,X_n)}
= \frac{V_{n,0}}{V_{n,n}} = \frac{V_{n,0}}{V_{n,m}}
\frac{V_{n,m}}{V_{n,n}}. \enas

Given $\epsilon \in (0,1)$, let $m$ be such that (\ref{Vnkinfsup}) holds.
But for any fixed $m$ we have by Theorem \ref{theorem-Ualpha} that
$$
\lim_{n \rightarrow \infty} \frac{V_{n,0}}{V_{n,m}}= 1.
$$

Hence by Lemma \ref{Vnn},
$$
\frac{1}{1+\epsilon} \le \liminf_{n \rightarrow \infty} \frac{V_{n,0}}{V_{n,n}}
 \le \limsup_{n \rightarrow \infty} \frac{V_{n,0}}{V_{n,n}} \le \frac{1}{1-\epsilon},
$$
and therefore the limit of the ratio exists and equals one.
Applying Theorem \ref{theorem-Ualpha} to the sequence
$n^{1/\alpha}V_{n,0}$ completes the proof of Theorem \ref{nancy}.
$\qed$

\section{Numerical Results and Additional Remarks}
%\label{table}
In Table 1, for the $\alpha = 0.1, 0.2, \ldots 1,
2, \ldots 10$ values in column (1), we tabulate the following
quantities in the columns indicated:\\

(2) $b_{\alpha}$

\vskip .5pc

(3) $\lim_{n \rightarrow \infty} nF (V^1_n) = (1 + 1/\alpha)$

\vskip .5pc

(4) $\lim_{n \rightarrow \infty} nF (V^2_n) =
h^{\alpha}(b_{\alpha}) = d^{\alpha}_{\alpha}$ and

\vskip .5pc

(5) $\lim_{n \rightarrow \infty} nF (\Vp)= \Gamma(1+1/\alpha)^{\alpha}$,\\

\noindent for $F(x)=x^\alpha L(x)$ and $\lim_{x \rightarrow 0}
L(x)=\B$ existing in $(0,\infty)$. In columns (6),(7), and (8), we
tablulate the ratios (3)/(4), (4)/(5) and (3)/(5).  Note that
another natural comparison would be among the values listed raised
to the power $1/\alpha$, as this would yield a comparison of the
actual limiting values of $V^1_n/V^2_n, V^2_n/E(\Vp)$ and
$V^1_n/E(\Vp)$ respectively. The reason that Table 1 lists the
values in the way it does is to display them in a comparable order
of magnitude to make numerical comparisons easier. The final
column of Table 1 presents the relative improvement attained by
using two stops rather than one, as compared to the reference
value of the prophet, \beq \label{eqno7.1} \lim_{n\to \infty}
(V^1_n-V^2_n)/(V^1_n-\Vp). \enq As evident from the table, the
improvement is highly significant for all values of $\alpha$.

The following asymptotic results can be shown to hold:

(i) For $\alpha \rightarrow \infty$,

\vskip  .5pc

$$\lim_{\alpha \rightarrow \infty} \lim_{n \rightarrow \infty} nF(V^1_n) = 1$$

\vskip .5pc

$$\lim_{\alpha \rightarrow \infty} \lim_{n \rightarrow \infty} nF(V^2_n) = 1-1/e$$

\vskip .5pc

$$\lim_{\alpha \rightarrow \infty}
\lim_{n \rightarrow \infty} nF(\Vp)= e^{-\gamma}$$ where
$\gamma=.5772 \ldots$ is Euler's constant. The limiting value for
the relative improvement (\ref{eqno7.1}) given in the last column
is

$$[1-\log(e-1)]/\gamma = 0.7946 \ldots $$

(ii) For $\alpha \to 0$,

The quantities in columns (3), (4) and (5) all tend to infinity, but
the ratios in columns (6),(7),(8) and (9) tend to a finite limit, and
are respectively
\vskip  .5pc

$$\lim_{\alpha \rightarrow 0} \lim_{n \rightarrow \infty}
\frac{nF(V^1_n)}{nF(V^2_n)}  = 2$$

\vskip .5pc

$$\lim_{\alpha \rightarrow 0} \lim_{n \rightarrow \infty}
\frac{nF(V^2_n)}{nF(\Vp)}  = e/2 = 1.3591  \ldots$$

\vskip .5pc

$$\lim_{\alpha \rightarrow 0} \lim_{n \rightarrow \infty}
\frac{nF(V^1_n)}{nF(\Vp)}  = e = 2.7182 \ldots$$

\vskip .5pc

The relative improvement (\ref{eqno7.1}) given in the last column
can be shown to tend to 1.

% comment out the \vskip line with a % and watch the table move.
% Use \begin{table}[h] to fix location
\vskip 2pc \oddsidemargin -.6in
\begin{table}
\begin{center}
\caption{Limiting Values of $nF(V^1_n), nF(V^2_n),
   nF(\Vp)$, and their ratios.}
\begin{tabular}{||c|c|c|c|c|c|c|c|c|c||}
\hline
(1) & (2) & (3) & (4) & (5) & (6) & (7) & (8) & (9)\\
\hline $\alpha$ & $b_{\alpha}$ & $\lim nF(V^1_n)$ & $\lim
nF(V^2_n)$ & $\lim nF(\Vp)$ &(3)/(4) &(4)/(5) &(3)/(5) & Eq. (\ref{eqno7.1}) \\
\hline
0.1 & 11.9312 & 11.0000 & 5.72334 & 4.52873 & 1.92195
& 1.26379 & 2.42894 &  .99868\\ \hline
0.2 & 6.8927 & 6.0000 & 3.20772 & 2.60517 & 1.87049
& 1.23129 & 2.30311 &  .97131 \\  \hline
0.3 & 5.2004 & 4.3333 & 2.36372 & 1.94980 & 1.83327
& 1.21229 & 2.22245 &  .93248 \\ \hline
0.4 & 4.3485 & 3.5000 & 1.93919 & 1.61670 & 1.80488
& 1.19947 & 2.16490 &  .90235 \\ \hline
0.5 & 3.8342 & 3.0000 & 1.68310 & 1.41421 & 1.78242
& 1.19013 & 2.12132 &  .88102 \\ \hline
0.6 & 3.4896 & 2.6667 & 1.51157 & 1.27776 & 1.76417
& 1.18298 & 2.08699 &  .86571\\ \hline
0.7 & 3.2423 & 2.4286 & 1.38853 & 1.17940 & 1.74902
& 1.17732 & 2.05916 &  .85460 \\ \hline
0.8 & 3.0561 & 2.2500 & 1.29590 & 1.10506 & 1.73624
& 1.17270 & 2.03610 &  .84614 \\ \hline
0.9 & 2.9107 & 2.1111 & 1.22362 & 1.04684 & 1.72530
& 1.16887 & 2.01665 &  .83958 \\ \hline
1.0 & 2.7940 & 2.0000 & 1.16562 & 1.00000 & 1.71583
& 1.16562 & 2.00000 &  .83438 \\ \hline
2.0 & 2.2634 & 1.5000 & 0.90214 & 0.78540 & 1.66270
& 1.14864 & 1.90984 &  .81217 \\ \hline
3.0 & 2.0839 & 1.3333 & 0.81309 & 0.71207 & 1.63983
& 1.14186 & 1.87245 & .80556 \\ \hline
4.0 & 1.9934 & 1.2500 & 0.76825 & 0.67497 & 1.62707
& 1.13820 & 1.85193 &   .80252 \\ \hline
5.0 & 1.9388 & 1.2000 & 0.74123 & 0.65255 & 1.61895
& 1.13590 & 1.83897 &  .80078 \\ \hline
6.0 & 1.9023 & 1.1666 & 0.72316 & 0.63753 & 1.61324
& 1.13432 & 1.82994 &  .79967 \\ \hline
7.0 & 1.8762 & 1.1429 & 0.71023 & 0.62677 & 1.60914
& 1.13317 & 1.82343 &  .79892 \\ \hline
8.0 & 1.8566 & 1.1250 & 0.70052 & 0.61867 & 1.60592
& 1.13230 & 1.81839 &  .79831 \\ \hline
9.0 & 1.8412 & 1.1112 & 0.69296 & 0.61236 & 1.60350
& 1.13162 & 1.81455 &  .79789 \\ \hline
10.0 & 1.8291 & 1.1000 & 0.68689 & 0.60731 & 1.60147
& 1.13105 & 1.81134 &  .79756\\ \hline \hline
\end{tabular}
\end{center}
\end{table}

\vskip 2pc
\begin{remark}
%\label{general.s.v}
Though we have proven Theorem \ref{nancy} for
the case where $F(x)=x^\alpha L(x)$, $\alpha>0$ and $L(x)$ having
finite positive limit as $x \downarrow 0$, we believe it holds
true for all $F \in {\cal D}(G^\alpha)$ of (\ref{Galpha}), that
is, whenever $L(x)$ is slowly varying as $x \downarrow 0$.
\end{remark}

\begin{remark}
\label{remark7.1} The approach in the present paper can easily be
applied to obtain the asymptotic behavior of the one-choice value
(obtained in Kennedy and Kertz (1991) by a different method), when
$F(x)=x^{\alpha}L(x)$ and $\lim_{x \downarrow 0} L(x)=\B\in(0,
\infty)$. First assume that $X \sim {\cal U}^{\alpha}(x)$ as in
(\ref{1.14}).  Then for the one choice value $V^1_n$, we have
$$
V^1_{n+1}=E[X \wedge V^1_n]=\alpha\int^{V^{1}_{n}}_0
x^{\alpha}dx+(1-(V^1_n)^{\alpha})V^1_n.
$$
Set $W^1_n=n^{1/\alpha}V^1_n$, and make the change of
variable $y=nx^{\alpha}$, as in Section
2. Now multiply by $n^{1/\alpha}$ to obtain
$$
\begin{array}{lll}
\left(\displaystyle\frac{n}{n+1}\right)^{1/\alpha} W^1_{n+1} &=&
\displaystyle\frac{1}{n} \int^{(W^1_n)^{\alpha}}_0 y^{1/\alpha} dy
+
(1-(V^1_n)^{\alpha}) W_n^1\\
\\
&=& \displaystyle\frac{1}{n}\int^n_0 (W^1_n \wedge
y^{1/\alpha}) dy.
\end{array}
$$
Thus $W^1_n$ satisfies (\ref{n4.5}) with $q(y)=y^{1/\alpha}$, and
now Theorem \ref{theorem4.1} can be applied to yield that $W^1_n
\to q(b)$ where $b$ is the unique root of
$$
\int^y_0 u^{1/\alpha} du + (1/\alpha-y)y^{1/\alpha}=0,
$$
giving $b=1+1/\alpha$. Hence, $\lim_{n\to \infty}
W^1_n=(1+1/\alpha)^{1/\alpha}$, or, $\lim_{n\to \infty}
nF(V^1_n)=(1+1/\alpha)$. The general result for the wider class of
distribution functions mentioned now follows in a manner similar
to, but simpler than, the calculation for two choices.
\end{remark}

\begin{remark}
%\label{remark7.2}
A similar approach can also be used to obtain
the limiting value for more than 2 choices. For three choices one
must first obtain the function $h^{(3)}(y)$ which replaces the
function $h^{(2)}(y)= h(y)$ of (8).  (Note that by Remark
\ref{remark7.1}, $h^{(1)}(y)=y^{1/\alpha})$.
\end{remark}

\begin{remark}
%\label{remark7.4}
Our results translate easily to the case where
the statistician is given two choices and his goal is to pick as
large a value as possible, his payoff being the expectation of the
larger of the two values chosen. Denote the optimal two-choice
value based on $n$ i.i.d. observations by ${\tilde V}^2_n$. Then
for $X \sim F(x)$, where $x_F < \infty$, and

$$
F_X(x) = 1 -(x_F - x)^\alpha L(x_F -x)
$$
where $L(\cdot)$ satisfies $lim_{y \downarrow 0} L(y) = \B$ and $0
< \B < \infty$, we have

$$
\lim_{n \rightarrow \infty}n[1 - F( {\tilde V}^2_n)] =
h^\alpha(b_\alpha).
$$
\end{remark}

\ignore{\begin{remark}
%\label{remark7.3}
Suppose $X$ belongs to the domain of attraction (for the maximum) of the
extremal distribution function $\exp(-e^{-x})$ for
$-\infty < x < \infty$.  A typical example is the
exponential distribution with mean 1, which we consider
below.  The aim is finding the optimal value of a
stopping rule which {\em maximizes} the value
when stopped.  As in Section 1, let $\hat{V}^1_n$ and
$\hat{V}^2_n$ denote the optimal values for 1 and 2
choices respectively, and $n$ observations, and let
$M_n = \max(x_1, \ldots, x_n)$.  Then it follows from
Kennedy and Kertz (1991) that for this case
$$
\lim_{n\to \infty}(EM_n - \log n)= \gamma =
\mbox{ Euler's constant} = .5772\ldots \eqno(7.4)
$$
and
$$
\lim_{n\to \infty}(\hat{V}^1_n - \log n) = 0.
\eqno(7.5)
$$
Thus we expect to have
$$
\lim_{n\to \infty}(\hat{V}^2_n - \log n) = c
\eqno(7.6)
$$
for some $0 < c < \gamma$.  The, maybe surprising fact
is that $c=0$, i.e. asymptotically one is no better off
with two choices than with one choice!  Clearly a
similar result holds for the other distributions in the
some domain of attraction.
\end{remark}}

\section{Final Remarks}
The last two authors are very saddened to announce that our
invaluable colleague and friend David Assaf passed away most
suddenly on December 23$^{rd}$ 2003 as this work was nearing
completion. On that very day, in a last email from Prof. Assaf to
us regarding the final touches on this manuscript, he wrote that
he had some ideas and `I will say more on this in a few days.' We
regret on many levels that this work can now only remain more or
less in its current form, without the benefit of those further
comments, now forever lost, which would have certainly greatly
improved the work.

\section*{Bibliography}
\begin{enumerate}

\item
%\label{Assaf1}
Assaf, D., and Samuel-Cahn, E. (2000)
Simple ratio prophet inequalities for a mortal with multiple choices.
{\em J. Appl. Prob.}, {\bf 37}, pp. 1084-1091.

\item
%\label{Assaf2}
Assaf, D., Goldstein, L. and Samuel-Cahn, E.  (2002).  Ratio
prophet inequalities when the mortal has several choices. {\em
Ann. Appl. Prob.}, {\bf 12}, pp. 972-984.

\item
%\label{de Haan}
de Haan, L.  (1976).  Sample extremes: an elementary introduction.
    {\em Statist.  Neerlandica}, {\bf 30}, pp. 161-172.

\item
%\label{KK90}
Kennedy, D.P. and Kertz, R.P. (1990). Limit theorems for threshold-stopped random
variables with applications to optimal stopping. {\em Adv. Appl.
Prob.} {\bf 22}, pp. 396-411.

\item
%\label{KK91}
Kennedy, D.P. and Kertz, R.P.  (1991).  The asymptotic behavior of
    the reward sequence in the optimal stopping of i.i.d random
    variables.  {\em Ann. Prob.} {\bf 19}, pp. 329-341.

\item
%\label{LLR}
Leadbetter, M.R., Lindgren, Georg and Rootz\'en, Holger.(1983).
{\em Extremes and Related Properties of Random Sequences and Processes}.
Springer-Verlag, New York, Heidelberg, Berlin

\item
%\label{Resnick}
Resnick, Sidney I. (1987).
{\it Extreme Values, Regular variation, and Point Processes.} Springer. N.Y

%\item \label{Peirce} Peirce, B.O. and Foster, R.M.  (1956).  A short table of
%integrals. Ginn and Company, Boston

\end{enumerate}
\end{document}